\documentclass[12pt,a4paper]{article}
\usepackage[T1]{fontenc}
\usepackage{lmodern}
\usepackage{amsmath,amssymb,amsthm,mathrsfs}
\usepackage{geometry,graphicx,booktabs}
\usepackage{pgfplots}
\pgfplotsset{compat=1.18}
\usepackage[hidelinks]{hyperref}
\hypersetup{pdftitle={The (p,q) norm of 2x2 circulant matrices II: non-Fourier extremals},pdfsubject={English version}}
\allowdisplaybreaks[2]
\newtheorem{theorem}{Theorem}[section]

\newtheorem{lemma}[theorem]{Lemma}

\theoremstyle{remark}
\newtheorem{remark}[theorem]{Remark}
\newtheorem*{analyticift}{Real analytic implicit function theorem}
\numberwithin{equation}{section}
\DeclareMathOperator{\arctanh}{arctanh}
\DeclareMathOperator{\arcosh}{arcosh}
\DeclareMathOperator{\arsinh}{arsinh}
\DeclareMathOperator{\sech}{sech}
\DeclareMathOperator{\sgn}{sgn}
\newcommand{\abs}[1]{\lvert#1\rvert}
\newcommand{\norm}[1]{\lVert#1\rVert}
\newcommand{\RR}{\mathscr R}
\newcommand{\vect}[1]{\boldsymbol{#1}}
\newcommand{\C}{\mathbb C}
\title{The $(p,q)$ Norm of $2\times2$ Circulant Matrices II\\[.3em]\large Non-Fourier Extremals}
\author{Di Wu\thanks{School of Science, Zhejiang University of Science and Technology,
Hangzhou 310023, P.R. China. E-mail address: \texttt{wudi08@mails.ucas.ac.cn}.
Supported by the Natural Science Foundation of Zhejiang Province of China (Grant
No.~LQ18A010002) and the National Natural Science Foundation of China (Grant
No.~12001488).}}
\date{}

\begin{document}
\maketitle
\begin{abstract}
We study the extremal structure of $2\times2$ complex circulant matrices whose induced $(p,q)$ norm is attained in neither Fourier direction. After normalization we restrict attention to parameters $\omega$ satisfying $0<|\omega|<1$, lying off both coordinate axes and outside the Weissler region for the two-point operator. For $1<p,q<2$ and for $2<p,q<\infty$ the open candidate sector contains exactly one stationary point, and that point is the global maximizer. For $1<p<2=q$ and for $p=2<q<\infty$ the problem reduces to a strictly monotone scalar equation. For $1<p<2<q<\infty$ the candidate sector contains at most three stationary points; when there are three, they are, in order of increasing radius, a strict local maximum, a saddle point and a strict local maximum, and whether $p<q'$ or $p>q'$ decides which peak is higher. For fixed $q$ and fixed matrix, a criterion at $p=q'$ tells whether some interval of $p$-values gives exactly three distinct stationary points in the open candidate sector; two turning values then separate the cases of one, two and three stationary points. We also characterize the norm implicitly at the dual exponent and give closed formulas for $p'=q=3$ and $4$.
\end{abstract}

\section{Introduction}

On the two-point probability space $\{-1,1\}$ let $d\nu:=(\delta_{-1}+\delta_1)/2$, and let $\varepsilon$ denote the coordinate function. For $a,b,\omega\in\C$ put $T_\omega(a+b\varepsilon):=a+\omega b\varepsilon$, and let $N_{p,q}(\omega):=\norm{T_\omega}_{L^p(d\nu)\to L^q(d\nu)}$ be the optimal constant. The two-point inequality reads
\begin{equation}\label{eq:two-point}
\left(\frac{|a+\omega b|^q+|a-\omega b|^q}{2}\right)^{1/q}
\le N_{p,q}(\omega)
\left(\frac{|a+b|^p+|a-b|^p}{2}\right)^{1/p}.
\end{equation}
The constant direction gives $N_{p,q}(\omega)\ge1$, so the inequality is contractive exactly when $N_{p,q}(\omega)=1$.

Let $\mathsf C$ be a nonzero $2\times2$ complex circulant matrix, with eigenvalues $\lambda_0$ and $\lambda_1$ on the constant and on the alternating Fourier direction. Conjugation by $\operatorname{diag}(1,-1)$ swaps the two directions and is an isometry on every $\ell^a$, so we may assume $|\lambda_0|\ge|\lambda_1|$ and put $\omega:=\lambda_1/\lambda_0$. In evaluation coordinates,
\[
M_\omega:=\begin{pmatrix}1+\omega&1-\omega\\1-\omega&1+\omega\end{pmatrix},
\qquad T_\omega=\tfrac12M_\omega,\qquad
\mathsf C=\tfrac{\lambda_0}{2}M_\omega .
\]
Since $\norm{f}_{L^a(d\nu)}=2^{-1/a}\norm{f}_{\ell^a}$,
\begin{equation}\label{eq:matrix-scaling}
\norm{\mathsf C}_{\ell^p\to\ell^q}
=|\lambda_0|2^{1/q-1/p}N_{p,q}(\omega).
\end{equation}
After normalization to the probability measure, the two Fourier directions give $1$ and $|\omega|$. Hence a Fourier extremal exists after normalization exactly when \eqref{eq:two-point} is contractive.

Two-point inequalities connect to the Hermite semigroup, complex hypercontractivity and optimal Fourier inequalities: Nelson's real-parameter hypercontractivity~\cite{Nelson1973} and Beckner's two-point method~\cite{Beckner1975} are the classical links. For $1<p\le q<\infty$, Weissler~\cite{Weissler1979} proved the result outside $3/2<p\le q<2$ and $2<p\le q<3$. Ivanisvili--Nazarov~\cite{IvanisviliNazarov2022} settled the remaining diagonal case, and a preprint of Huang--Ivanisvili~\cite{HuangIvanisvili2026} supplies the two-point contraction theorem in the strictly off-diagonal case. The first paper of this series~\cite{PartI} gives a different proof via ridge and coupling lines and also covers $p>q$.

Below, $a':=a/(a-1)$ denotes the conjugate exponent and $\overline{\mathbb D}:=\{z:|z|\le1\}$ is the closed unit disc. For $1<p,q<\infty$ put $\Omega_{p,q}:=\{\omega:N_{p,q}(\omega)=1\}$ and
\begin{equation}\label{eq:weissler}
A_{p,q}:=\{\omega:\ |(p-2)-(q-2)\omega^2|\le p-q|\omega|^2\}.
\end{equation}
\begin{theorem}[Fourier extremal criterion from the first paper]\label{thm:prior}
Let $1<p,q<\infty$. Then $\Omega_{p,q}=A_{p,q}\cap\overline{\mathbb D}$. In particular, if $|\lambda_0|\ge|\lambda_1|$, the matrix $\mathsf C$ admits a Fourier extremal if and only if $\lambda_1/\lambda_0\in A_{p,q}$.
\end{theorem}

On nonzero inputs we use the projective coordinate $\vect f=(1,z)^\top$, allowing $z=\infty$ when the first component vanishes. In the non-axis range considered below, and with the normalization of~\cite{PartI}, set
\begin{equation}\label{eq:ratio}
c:=\frac{1-\omega}{1+\omega},\qquad
\RR(z):=\frac{2^{1/p-1/q}}{\abs{1+c}}\,
\frac{(\abs{1+cz}^q+\abs{c+z}^q)^{1/q}}{(1+\abs z^p)^{1/p}}.
\end{equation}
Here $\RR(z)=\norm{T_\omega(1,z)^\top}_q/\norm{(1,z)^\top}_p$, since $1+c=2/(1+\omega)$. Once $\norm f_p=1$, the inputs lie in a compact set, so the maximum is attained. Since $\RR(z)=\RR(1/z)$, a maximizer can always be chosen in the closed unit disc. The two Fourier directions give $\RR(1)=1$ and $\RR(-1)=\abs\omega$.

Throughout, $1<p,q<\infty$, and we stay away from the axis and from the unitary boundary. Modulation $\omega\mapsto-\omega$ and complex conjugation both preserve the norm and produce the corresponding input, so we may assume
\begin{equation}\label{eq:domain}
\omega=s e^{i\phi},\quad 0<s<1,\quad 0<\phi<\frac{\pi}{2},
\qquad c=\rho e^{-i\gamma},\quad 0<\rho<1,\quad 0<\gamma<\frac{\pi}{2}.
\end{equation}
The ranges for $c$ follow from $\abs{1-\omega}<\abs{1+\omega}$ together with
\[
\operatorname{Re}c=\frac{1-s^2}{\abs{1+\omega}^2},
\qquad
\operatorname{Im}c=-\frac{2\operatorname{Im}\omega}{\abs{1+\omega}^2}.
\]
We say that $\omega\notin A_{p,q}$ lies outside the Weissler region. To see what this local condition means, put $w:=(1-z)/(1+z)=w_1+iw_2$; the two-point input is proportional to $1+w\varepsilon$. From
\[
\frac{1}{a}\log\frac{\abs{1+w}^a+\abs{1-w}^a}{2}
=\frac{(a-1)w_1^2+w_2^2}{2}+O(\abs w^4)
\]
we get the expansion around the constant direction. Indeed, writing $\abs{1\pm w}^a=(1\pm2w_1+\abs w^2)^{a/2}$, the odd-order terms cancel after averaging, the second-order term is $a\abs w^2/2+a(a-2)w_1^2/2$, and taking the logarithm and dividing by $a$ gives the display. With $a_0:=\operatorname{Re}\omega$ and $b_0:=\operatorname{Im}\omega$ this yields
\begin{equation}\label{eq:local}
\begin{gathered}
\log\RR(z)=\frac12(Ew_1^2+2Kw_1w_2+Dw_2^2)+O(\abs w^4),\\
E:=s^2+(q-2)a_0^2-(p-1),\quad
K:=-(q-2)a_0b_0,\quad
D:=s^2+(q-2)b_0^2-1.
\end{gathered}
\end{equation}
At $(r,\theta)=(1,0)$ the real differential of $w=(1-re^{i\theta})/(1+re^{i\theta})$ is $-\tfrac12 I$. Since $\RR(1)=1$ and the first differential vanishes, $\log\RR$ and $\RR$ have the same Hessian there. In polar coordinates $(r,\theta)$,
\[
\operatorname{Hess}_{(r,\theta)}\RR(1,0)
=\frac14\begin{pmatrix}E&K\\K&D\end{pmatrix}.
\]
In particular $D/4$ is the second angular derivative at $\theta=0$ on the unit circle. The trace of the matrix $\left(\begin{smallmatrix}E&K\\K&D\end{smallmatrix}\right)$ is $qs^2-p$, and
\[
E-D+i\,2K=(q-2)\overline{\omega}^{\,2}-(p-2).
\]
Its larger eigenvalue is $\bigl(qs^2-p+\abs{(p-2)-(q-2)\omega^2}\bigr)/2$, so negative semidefiniteness is equivalent to \eqref{eq:weissler}. Outside the region the quadratic form in \eqref{eq:local} is strictly positive in some direction, and a small nonzero $w$ in that direction gives $N_{p,q}(\omega)>1$. Theorem~\ref{thm:prior} then says that, under \eqref{eq:domain}, the outside case is precisely the one in which neither Fourier direction attains the norm. The three theorems below are stated under \eqref{eq:domain} and outside the Weissler region.

\begin{theorem}[$1<p,q<2$ or $2<p,q<\infty$]\label{thm:same}
If $1<p,q<2$, then $0<r<1$, $-\gamma<\theta<0$ contains exactly one two-dimensional stationary point. If $2<p,q<\infty$, then $0<r<1$, $0<\theta<\gamma$ contains exactly one two-dimensional stationary point. In both cases this stationary point is the unique global maximizer in the closed unit disc.
\end{theorem}

\begin{theorem}[$1<p<2=q$ or $p=2<q<\infty$]\label{thm:edges}
If $1<p<2=q$, the unique global maximizer in the closed unit disc lies on the positive real axis, and its radius is the unique root of \eqref{eq:edge-root} in $0<r<1$. If $p=2<q<\infty$, the unique global maximizer in the closed unit disc is recovered from \eqref{eq:edge-recover}.
\end{theorem}

\begin{theorem}[$1<p<2<q<\infty$]\label{thm:main}
In $0<r<1$, $0<\theta<\gamma$ there are one, two or three distinct stationary points.
\begin{enumerate}
\item If there is one, it is the global maximizer. If there are two, one is a strict local maximum and is the global maximizer, and the other is a horizontal inflection point along the corresponding profile.
\item If there are three, they are, in order of increasing radius, a strict local maximum, a saddle point and a strict local maximum. The inner peak is the higher one when $p<q'$, the outer peak when $p>q'$, and the two peaks are of equal height when $p=q'$.
\item Fix $q$ and $\omega$. Among the $p$-values for which $\omega\notin A_{p,q}$, there is an interval with exactly three distinct stationary points in the open candidate sector if and only if there is a positive diagonal point $d$ satisfying \eqref{eq:threshold}. If such an interval exists, the values of $\mathcal P$ at its two turning points determine two critical exponents $p_-<p_+$. For $p$ with $\omega\notin A_{p,q}$, there are three stationary points in the open interval $(p_-,p_+)$, two at either endpoint, and one outside the closed interval. If no such interval exists, there is exactly one stationary point for every $p$ with $\omega\notin A_{p,q}$.
\end{enumerate}
\end{theorem}

In each case, a global maximizer outside the disc is obtained from the one inside by $z\mapsto1/z$.

\paragraph{Method and structure.}
Inside the region where maximizers may lie, fix the radius of the input and maximize in the angle. The resulting curve, the ridge line, depends only on the output exponent $q$. Substituting the weight ratio from the angular stationarity equation into the radial equation gives the coupling line, which depends only on the input exponent $p$. The two-dimensional stationary points in the candidate sector are the intersections of these two curves.

Section~\ref{sec:prelim} establishes uniqueness of the ridge line, the parametric description of the coupling line and the location of its exit point. Theorem~\ref{thm:same} is proved in Section~\ref{sec:sameside}. For $1<p,q<2$ uniqueness along the coupling line follows from the intersection sign criterion of the first paper; for $2<p,q<\infty$ the dual correspondence of stationary points reduces the claim to the previous case. Theorem~\ref{thm:edges} is proved in Section~\ref{sec:edges}, where $1<p<2=q$ reduces to a monotone equation for the radius and $p=2<q<\infty$ follows by duality. Theorem~\ref{thm:main} is proved in Section~\ref{sec:count}. There the same pair of curves is written in weighted coordinates, and the number of intersections becomes the number of intersections of a function of one variable with a horizontal line. Section~\ref{app:closed} gives closed formulas for the norms at $(4/3,4)$ and $(3/2,3)$.

\section{Preliminaries}\label{sec:prelim}

\subsection{Candidate sector, ridge line and coupling equation}\label{sec:polar}
Write $z=re^{i\theta}$ and put
\begin{equation}\label{eq:UV}
U:=1+\rho^2r^2+2\rho r\cos(\gamma-\theta),\qquad
V:=\rho^2+r^2+2\rho r\cos(\gamma+\theta).
\end{equation}
On a fixed circle, maximizing $\RR$ is the same as maximizing $S_r(\theta):=U^{q/2}+V^{q/2}$. Where $U,V>0$,
\begin{equation}\label{eq:ridge}
S_r'(\theta)=q\rho r
\bigl[U^{(q-2)/2}\sin(\gamma-\theta)-V^{(q-2)/2}\sin(\gamma+\theta)\bigr].
\end{equation}
If one of the squared moduli vanishes, the corresponding $q$-th power still has zero first derivative at that point, because $q>1$. Hence the angular derivative is nonpositive on $[\gamma,\pi-\gamma]$ and nonnegative on $[-\pi+\gamma,-\gamma]$, and strictly so at interior points of either interval. Write the left sector as $\theta=\pi+\psi$ with $\abs\psi\le\gamma$. If $\cos(\psi-\gamma)$ and $\cos(\psi+\gamma)$ are both nonnegative, put $\theta_0=\psi$; then neither squared modulus decreases and at least one of them strictly increases. If $\cos(\psi-\gamma)<0$, put $\theta_0=2\gamma-\pi-\psi$; then $\theta_0\in[-\gamma,\gamma]$ and
\[
U(r,\theta_0)=U(r,\pi+\psi),\quad
V(r,\theta_0)-V(r,\pi+\psi)=-4\rho r\sin(2\gamma)\sin(\psi-\gamma)>0.
\]
If $\cos(\psi+\gamma)<0$, put $\theta_0=\pi-2\gamma-\psi$, and again $V$ is unchanged while $U$ strictly increases. The two cosines cannot both be negative, since that would require $\psi<\gamma-\pi/2$ and $\psi>\pi/2-\gamma$. This confines the per-circle maximum to $[-\gamma,\gamma]$.

For $0<r<1$ one has $1+\rho^2r^2>\rho^2+r^2$, and for $0<\theta<\gamma$ one has $\cos(\gamma-\theta)>\cos(\gamma+\theta)$. Write
\[
m_1:=1+\rho^2r^2,\quad m_2:=\rho^2+r^2,\quad
d_+:=2\rho r\cos(\gamma-\theta),\quad d_-:=2\rho r\cos(\gamma+\theta),
\]
so that $m_1>m_2$ and $d_+>d_-$, while $(U,V)(r,\theta)=(m_1+d_+,m_2+d_-)$ and $(U,V)(r,-\theta)=(m_1+d_-,m_2+d_+)$. The difference between $S_r(\theta)$ and $S_r(-\theta)$ therefore lies only in whether the larger constant term is paired with the larger cosine.

\begin{lemma}[Location of the per-circle maximum]\label{lem:sector-compare}
Let $0<r<1$, $q>1$, $q\ne2$, and write $f(t):=t^{q/2}$. Then for $0<\theta<\gamma$,
\[
S_r(\theta)-S_r(-\theta)
=\bigl[f(m_1+d_+)-f(m_1+d_-)\bigr]-\bigl[f(m_2+d_+)-f(m_2+d_-)\bigr],
\]
and this is strictly positive for $q>2$ and strictly negative for $q<2$. Hence the candidate sector is $-\gamma<\theta<0$ for $q<2$ and $0<\theta<\gamma$ for $q>2$, in both cases with $0<r<1$.
\end{lemma}
\begin{proof}
For $q>2$ the function $f$ is strictly convex, so the difference quotient $\bigl(f(x+h)-f(x)\bigr)/h$ strictly increases with the base point $x$. Take the two base points $m_2+d_-<m_1+d_-$ with the common increment $d_+-d_->0$; then $f(m_1+d_+)-f(m_1+d_-)>f(m_2+d_+)-f(m_2+d_-)$, that is, the difference is strictly positive. For $q<2$ the function $f$ is strictly concave and the inequality reverses. For $q=2$, $U+V$ depends on the angle only through a positive coefficient times $\cos\theta$, the per-circle maximum is at $\theta=0$, and this case is treated separately.
\end{proof}
From \eqref{eq:ridge} one also reads off $S_r'(-\gamma)>0$ and $S_r'(\gamma)<0$, so the derivatives at the two endpoints point into the sector; at $\theta=0$ the sign of $S_r'(0)$ is that of $q-2$. Together with Lemma~\ref{lem:sector-compare} this places the per-circle maximizer in the interior of the open sector.

\begin{lemma}[Location of the maximum]\label{lem:argmax}
Let $q\ne2$ and let $\omega$ lie outside the Weissler region. Then the maximum of $\RR$ over the closed unit disc is attained in the open candidate sector: $-\gamma<\theta<0$ for $q<2$ and $0<\theta<\gamma$ for $q>2$, in both cases with $0<r<1$.
\end{lemma}
\begin{proof}
\textit{Step 1: the circle.} If $q<2$, strict concavity of $t^{q/2}$ gives
\[
S_1(\theta)\le 2\Bigl[\frac{U+V}{2}\Bigr]^{q/2},
\qquad U+V=2(1+\rho^2)+4\rho\cos\gamma\cos\theta,
\]
and the right-hand side is uniquely maximal at $\theta=0$, where $U=V$. Thus $z=1$ is the unique maximizer on the circle, with value $\RR(1)=1$; since the maximum outside the region is strictly larger than $1$, no point of the circle is a maximizer. If $q>2$, use $\RR(z)=\RR(1/z)$ first: on the circle $1/z=\bar z$, so $\RR$ is invariant under $\theta\mapsto-\theta$ and it suffices to look at $\theta\ge0$. At $r=1$, differentiate
\[
\log\RR=\mathrm{const}+\frac1q\log(U^{q/2}+V^{q/2})-\frac1p\log(1+r^p)
\]
in the radial direction. Using $U_r-U=\rho^2-1$ and $V_r-V=1-\rho^2$ gives
\begin{equation}\label{eq:radial-circle}
\left.\partial_r\log\RR\right|_{r=1}
=\frac{(1-\rho^2)\big(V^{(q-2)/2}-U^{(q-2)/2}\big)}
{2(U^{q/2}+V^{q/2})}.
\end{equation}
Here $U-V=4\rho\sin\gamma\sin\theta$, so $U>V$ for $0<\theta<\pi$ and hence $V^{(q-2)/2}<U^{(q-2)/2}$; thus \eqref{eq:radial-circle} is strictly negative: moving inward along the radius strictly increases $\RR$, and the point of the circle is not a maximizer. At $\theta=0,\pi$ the values are $\RR(1)=1$ and $\RR(-1)=s<1$, both smaller than the maximum outside the region. So for $q>2$ the maximizer also lies in the open disc.

\textit{Step 2: inside the open disc.} A maximizer in the open disc is a two-dimensional stationary point, and its polar angle must maximize $S_r$ on the same circle. By \eqref{eq:ridge},
\[
S_r'(0)=q\rho r\sin\gamma\bigl(U^{(q-2)/2}-V^{(q-2)/2}\bigr),
\]
while $U-V=(1-r^2)(1-\rho^2)>0$, so $S_r'(0)\ne0$ for $0<r<1$ and $q\ne2$. Lemma~\ref{lem:sector-compare} and the signs of the derivatives at the two ends of the sector then put the polar angle of the maximizer in the corresponding open sector. Finally, $z=0$ is not a maximizer: along $z=t\ge0$ the first variation of the denominator vanishes, while that of the power sum in the numerator equals $q(1+\rho^{q-2})\operatorname{Re}c>0$. With the reciprocal symmetry $\RR(z)=\RR(1/z)$ this proves the assertion.
\end{proof}

Inside the open sector, setting the bracket in \eqref{eq:ridge} equal to zero gives the angular stationarity equation. The radial equation is
\[
\frac{U^{(q-2)/2}(\rho^2r+\rho\cos(\gamma-\theta))
+V^{(q-2)/2}(r+\rho\cos(\gamma+\theta))}
{U^{q/2}+V^{q/2}}=\frac{r^{p-1}}{1+r^p}.
\]
The angular equation fixes the ratio of the two positive weights to be $\sin(\gamma+\theta):\sin(\gamma-\theta)$. Substituting these two sines into the radial equation and using $\sin a\cos b+\cos a\sin b=\sin(a+b)$ turns it into
\begin{equation}\label{eq:coupling-polar}
\frac{r[\rho^2\sin(\gamma+\theta)+\sin(\gamma-\theta)]+\rho\sin(2\gamma)}
{(1+\rho^2 r^2)\sin(\gamma+\theta)+(\rho^2+r^2)\sin(\gamma-\theta)+2\rho r\sin(2\gamma)}
=\frac{r^{p-1}}{1+r^p}.
\end{equation}
Every quantity canceled in this step is positive, so the step is reversible. Equation \eqref{eq:coupling-polar} is called the coupling equation; it involves $p$ only. The ridge line formed by the unique per-circle maximizers involves $q$ only.

Put
\[
\begin{gathered}
\tau:=-\frac p2\log r,\quad \alpha_p:=\frac{2-p}{p},\quad\alpha_q:=\frac{2-q}{q},\\
\sigma:=\frac{2\rho}{1+\rho^2},\quad \delta:=\frac{1-\rho^2}{1+\rho^2},
\quad C:=\cos\gamma,\quad S:=\sin\gamma.
\end{gathered}
\]
The coupling equation is then equivalent to
\begin{equation}\label{eq:coupling-hyp}
S\cos\theta\,\sinh(\alpha_p\tau)
+\delta C\sin\theta\,\cosh(\alpha_p\tau)-\sigma SC\sinh\tau=0.
\end{equation}
To check the substitution, note that \(\rho^2\sin(\gamma+\theta)+\sin(\gamma-\theta)
=(1+\rho^2)(S\cos\theta-\delta C\sin\theta)\), and that the constant and the $r^2$ terms of the denominator equal $(1+\rho^2)(S\cos\theta+\delta C\sin\theta)$ and $(1+\rho^2)r^2(S\cos\theta-\delta C\sin\theta)$ respectively. After cross multiplication the common $r^{p+1}$ terms cancel and what remains is
\[
S\cos\theta(r-r^{p-1})-\delta C\sin\theta(r+r^{p-1})+\sigma SC(1-r^p)=0.
\]
Division by $-2r^{p/2}$ gives \eqref{eq:coupling-hyp}.

\subsection{The fourth-quadrant ridge line and the case \texorpdfstring{$q=2$}{q=2}}

For $q<2$ the stationarity equation \eqref{eq:ridge} can be recast as the angular residual used in Section~\ref{sec:sameside}. With $\xi=-\log r>0$, put
\[
X:=\frac{\delta\sinh\xi+\sigma S\sin\theta}{\cosh\xi+\sigma C\cos\theta}.
\]
Since $U/V=(1+X)/(1-X)$, one has $|X|<1$ in the open sector, and
\[
\frac{\partial X}{\partial\theta}
=\frac{\sigma[S\cos\theta\cosh\xi+\delta C\sin\theta\sinh\xi+\sigma SC]}
{(\cosh\xi+\sigma C\cos\theta)^2}>0 .
\]
Here $C|\sin\theta|<S\cos\theta$, so the bracket exceeds
$S\cos\theta(\cosh\xi-\delta\sinh\xi)>0$. Thus
\[
\alpha_q\arctanh X+(1+\alpha_q)\arctanh(C\tan\theta/S)
\]
strictly increases on $-\gamma<\theta<0$, from $-\infty$ to a positive value. The first limit comes from the second $\arctanh$ and the second from $X(\xi,0)>0$. The function has a unique zero and its derivative is strictly positive, so the ridge line is real analytic for $0<r<1$.

When $q=2$, the per-circle comparison of Section~\ref{sec:polar} gives $\theta_R(r)=0$, which is real analytic in $r$.

\subsection{The first-quadrant ridge line and the maximal angle on the unit circle}\label{sec:upper-ridge}
Let $q>2$ and fix $0<r<1$. On $0<\theta<\gamma$ write
\[
L(\theta):=\log\frac{\sin(\gamma+\theta)}{\sin(\gamma-\theta)},\qquad
L_V(\theta):=\log\frac{V(r,-\theta)}{V(r,\theta)}.
\]
By \eqref{eq:UV},
\[
\log\frac{U(r,\theta)}{V(r,\theta)}
=L_V(\theta)+\log\left(1+\frac{(1-\rho^2)(1-r^2)}{V(r,-\theta)}\right).
\]
The last term is strictly positive and strictly decreasing in $\theta$, because $V(r,-\theta)=\rho^2+r^2+2\rho r\cos(\gamma-\theta)$ is strictly increasing. After dividing both terms by $L$ we look at the first one. Direct differentiation gives
\[
L'(\theta)=\frac{\sin(2\gamma)}{\sin(\gamma-\theta)\sin(\gamma+\theta)}>0,\qquad
L_V'(\theta)=2\rho r\left[
\frac{\sin(\gamma-\theta)}{V(r,-\theta)}
+\frac{\sin(\gamma+\theta)}{V(r,\theta)}\right].
\]
On a common denominator, and using
\[
\sin(\gamma+\theta)\sin(\gamma-\theta)=\cos^2\theta-C^2,\quad
\cos(\gamma+\theta)\cos(\gamma-\theta)=\cos^2\theta+C^2-1
\]
together with the addition formula for sine, one obtains in terms of $u_\theta:=\cos\theta$
\[
\frac{L_V'}{L'}=
\frac{2\rho r}{C}\big[(\rho^2+r^2)u_\theta+2\rho r C\big]
\frac{u_\theta^2-C^2}{(\rho^2+r^2)^2+4\rho r(\rho^2+r^2)C u_\theta+4\rho^2r^2(u_\theta^2+C^2-1)}.
\]
The reciprocal of the second factor equals
\[
4\rho^2r^2+\frac{4\rho r(\rho^2+r^2)C u_\theta+(\rho^2+r^2)^2+4\rho^2r^2(2C^2-1)}{u_\theta^2-C^2}.
\]
The coefficient of the linear term in the numerator is positive, and the constant term is also strictly positive because $(\rho^2+r^2)^2\ge4\rho^2r^2$. On $u_\theta>C$ one has $u_\theta/(u_\theta^2-C^2)=1/(u_\theta-C^2/u_\theta)$. The denominator there and $u_\theta^2-C^2$ are both positive and strictly increasing, so $u_\theta/(u_\theta^2-C^2)$ and $1/(u_\theta^2-C^2)$ are both strictly decreasing. Hence both positive factors above strictly increase with $u_\theta$, that is, $L_V'/L'$ strictly decreases in $\theta$. Since $L_V(0)=L(0)=0$, Cauchy's mean value theorem gives, for every $0<\theta<\gamma$, a point $\widehat\theta\in(0,\theta)$ with
\[
\frac{L_V(\theta)}{L(\theta)}
=\frac{L_V'(\widehat\theta)}{L'(\widehat\theta)}
>\frac{L_V'(\theta)}{L'(\theta)}.
\]
Differentiating $L_V/L$ directly shows that the derivative equals $(L'/L)(L_V'/L'-L_V/L)$, hence is strictly negative. The other positive term divided by $L$ is also strictly decreasing, so the whole quotient $\log(U/V)/L$ strictly decreases. From the positivity of the logarithm as $\theta\downarrow0$ together with $L(0)=0$, and from the finiteness of the logarithm as $\theta\uparrow\gamma$ together with $L\to\infty$, the two limits are $+\infty$ and $0$.

The angular equation can be written as $\log(U/V)/L=2/(q-2)$, so it has a unique root. The angular residual is
\[
\frac{q-2}{2}\log(U/V)-L .
\]
At the root its angular derivative equals $\frac{q-2}{2}L$ times the strictly negative derivative above. Expression \eqref{eq:ridge} equals a positive factor times the difference between the exponential of the residual and $1$, so the second angular derivative at the root is strictly negative. The real analytic implicit function theorem gives real analyticity of the ridge line $\theta_R(r)$ for $0<r<1$. It is determined by $q$ and $\omega$.

At $r=1$ the positive correction term disappears. The limit of $\log(U/V)/L$ as $\theta=0^+$ is then given by $L_V'/L'$ and equals
\[
\frac{2\rho\sin^2\gamma}{C(1+\rho^2+2\rho C)}.
\]
Hence, if $(q-2)\rho\sin^2\gamma>C(1+\rho^2+2\rho C)$, the angular equation has a unique positive root in $0<\theta<\gamma$, and $S_1$ first rises and then falls on that interval. The two points symmetric about the real axis are the maximizers on the circle. If the inequality is reversed or equality holds, strict monotonicity forces the residual to be strictly negative on $0<\theta<\gamma$. The restriction on candidate angles in Section~\ref{sec:polar} and symmetry then give zero as the only maximizing angle. Using $1+c=2/(1+\omega)$ together with $\rho C=\operatorname{Re}c$ and $\rho\sin\gamma=-\operatorname{Im}c$, the strict inequality is equivalent to $(q-2)(\operatorname{Im}\omega)^2>1-s^2$, that is, to $D>0$ in \eqref{eq:local}. This also proves uniqueness when $D=0$.

\subsection{Existence, exit and radial folding of the coupling line}\label{app:coupling}

The coupling equation \eqref{eq:coupling-hyp} depends only on $p$ and $\omega$. Here we ask whether the radius can be used as a parameter; the relevant plots are in Appendix~\ref{app:coupling-plots}. The turning points of Section~\ref{sec:count} are the points of the weight plane where $\mathcal P'(x)=0$.

\paragraph{Angle as parameter.}

Let $\mathcal H(\tau,\theta)$ denote the left-hand side of \eqref{eq:coupling-hyp}. If $1<p<2$, fix $0<\theta\le\gamma$ and divide $\mathcal H$ by $\sinh\tau>0$. In the first term, $\sinh(\alpha_p\tau)/\sinh\tau$ strictly decreases, its logarithmic derivative being $\alpha_p\coth(\alpha_p\tau)-\coth\tau<0$ because $t\coth t$ strictly increases; the derivative of the latter is $(\sinh t\cosh t-t)/\sinh^2t>0$. In the second term, $\cosh(\alpha_p\tau)/\sinh\tau$ also strictly decreases, its logarithmic derivative being $\alpha_p\tanh(\alpha_p\tau)-\coth\tau<0$. Both coefficients are positive. After division by $\sinh\tau$ the sum of the first two terms decreases strictly from infinity to zero. Hence $\mathcal H=0$ has a unique positive root $\tau_C(\theta)$. At the root the product rule gives
\[
\mathcal H_\tau=\sinh\tau\ \partial_\tau(\mathcal H/\sinh\tau)<0,
\]
because the term multiplied by $\cosh\tau$ vanishes exactly. By the real analytic implicit function theorem, $r_C(\theta)=e^{-2\tau_C(\theta)/p}$ is single-valued and real analytic on $0<\theta\le\gamma$, and can be extended locally beyond $\gamma$ on either side in the angular variable.

If $p\ge2$, the same argument gives a unique root at fixed angle. If $p=2$, the first term in \eqref{eq:coupling-hyp} vanishes, the second does not depend on $\tau$, and the last term strictly decreases. If $p>2$ then $\alpha_p<0$, and the derivative in $\tau$ of the sum of the first two terms equals
\[
|\alpha_p|\big[-S\cos\theta\,\cosh(|\alpha_p|\tau)
+\delta C\sin\theta\,\sinh(|\alpha_p|\tau)\big]<0,
\]
because $\delta C\sin\theta/(S\cos\theta)\le\delta<1$. Adding the strictly negative derivative of the last term gives $\mathcal H_\tau<0$. Also $\mathcal H(0,\theta)=\delta C\sin\theta>0$, while $\mathcal H$ tends to $-\infty$ at the far end since $|\alpha_p|<1$, so the unique root exists.

\paragraph{The two boundary radii.}

Since $\sigma=1/\cosh(\log(1/\rho))$ and $\delta=\tanh(\log(1/\rho))$, the exit of the ray can be simplified with the logarithmic radius. Substituting $\theta=\gamma$, the coupling equation becomes
\[
SC\bigl[\sinh(\alpha_p\tau)+\delta\cosh(\alpha_p\tau)-\sigma\sinh\tau\bigr]=0.
\]
Multiplying the first two terms by $\cosh(\log(1/\rho))$ and using the addition formula for $\sinh$ gives
$\sinh(\alpha_p\tau+\log(1/\rho))=\sinh\tau$. Since $\sinh$ is strictly increasing on the real axis,
\[
\tau_C(\gamma)=\frac{\log(1/\rho)}{1-\alpha_p},\qquad
r_C(\gamma)=e^{-2\tau_C(\gamma)/p}=\rho^{1/(p-1)}.
\]

The boundary on the positive real axis is described by $\theta\downarrow0$. For $p<2$ the roots $\tau_C(\theta)$ are uniformly bounded: as $\tau\to\infty$ both positive terms of $\mathcal H/\sinh\tau$ tend to $0$ uniformly in $\theta$, so $\mathcal H/\sinh\tau\to-\sigma SC<0$. A positive subsequential limit would have to satisfy $\sinh(\alpha_p\tau)/\sinh\tau=\sigma C$. For $\alpha_p>\sigma C$ this equation has a unique positive root $\tau_*$, and the limit cannot be zero, because the normalized first term would tend to $S\alpha_p>\sigma SC$ while the second term is nonnegative. Hence $\tau_C(0^+)=\tau_*$. For $\alpha_p\le\sigma C$ the equation has no positive root, so $\tau_C(0^+)=0$. For $p\ge2$ one has $\mathcal H(\tau,0)<0$ for all $\tau>0$, and the roots at positive angles are uniformly bounded, so again $\tau_C(0^+)=0$. The boundary radius on the real axis is therefore either the exit radius $e^{-2\tau_*/p}$ of the fourth-quadrant branch or it is $1$.

\paragraph{One or two branches at fixed radius.}

For $p\ge2$ one has $\mathcal H_\theta>0$ on $0<\theta<\gamma$, because
\[
\mathcal H_\theta=-S\sin\theta\,\sinh(\alpha_p\tau)+\delta C\cos\theta\,\cosh(\alpha_p\tau).
\]
Hence $\tau_C'=-\mathcal H_\theta/\mathcal H_\tau>0$, that is, $r_C$ strictly decreases in $\theta$. The coupling line can therefore be written uniquely as a branch with $r$ as the independent variable, and it leaves the sector through $\theta=\gamma$.

For $1<p<2$, on the side $\theta\downarrow0$ the derivative $\mathcal H_\theta$ tends to $\delta C\cosh(\alpha_p\tau_C(0^+))>0$. If $\tau_C'(\theta)=0$, differentiating $\mathcal H(\tau_C(\theta),\theta)=0$ once more gives $\tau_C''=-\mathcal H_{\theta\theta}/\mathcal H_\tau$, since all other second-order chain rule terms contain $\tau_C'$. On the curve, $\mathcal H=0$ gives
\[
\mathcal H_{\theta\theta}
=-S\cos\theta\,\sinh(\alpha_p\tau)-\delta C\sin\theta\,\cosh(\alpha_p\tau)
=-\sigma SC\sinh\tau<0 .
\]
Thus $\tau_C''<0$, and $r_C''=-(2/p)r_C\tau_C''>0$. Every interior radial turning point is therefore a strict local minimum of the radius, and there is at most one in the interior.

Moreover, the sign of $\mathcal H_\theta$ at $\theta=\gamma$ equals the sign of
$\delta C^2-S^2\tanh(\alpha_p\tau_C(\gamma))$. If this is negative, the positive sign at the starting angle shows that there is exactly one interior turning point. If it is positive, a turning point would force the sign to its right to stay negative, a contradiction. If it vanishes, then at the endpoint of the ray $\tau_C'=0$ and $\tau_C''<0$, so a left neighborhood of that endpoint is still strictly increasing and there is no interior turning point either. The necessary and sufficient condition for two radial branches is therefore
\begin{equation}\label{eq:coupling-fold-criterion}
\tanh\left(\frac{\alpha_p(-\log\rho)}{1-\alpha_p}\right)
>\delta\,\frac{\cos^2\gamma}{\sin^2\gamma}.
\end{equation}
Under this condition, between the minimal radius and the smaller of the two boundary radii every radius corresponds to two angles: one branch runs to the real axis and the other to $\theta=\gamma$. Beyond the smaller of the two boundary radii, whether the circle of that radius meets the line once more is decided by the relative size of the two boundary radii. When there is no interior turning point, the radius on the whole branch decreases strictly with the angle, and its minimum is attained on the boundary ray. The two branches appear only when the radius is used as the parameter; parametrized by the angle there is always a single continuous curve.
\paragraph{The fourth quadrant: one branch and its exit.}
Let $1<p<2$. On $-\gamma<\theta<0$ one has $\mathcal H_\theta>0$, and at $\theta=-\gamma$
\[
\mathcal H(\tau,-\gamma)=SC[\sinh(\alpha_p\tau)-\delta\cosh(\alpha_p\tau)-\sigma\sinh\tau]<0 .
\]
Here one uses an identity: from $\delta^2+\sigma^2=1$ and the subtraction formula for $\sinh$,
\begin{equation}\label{eq:neg-bracket}
\sinh(\alpha_p\tau)-\delta\cosh(\alpha_p\tau)=\sigma\sinh(\alpha_p\tau-u),
\qquad u:=\log\frac1\rho=\arctanh\delta ,
\end{equation}
so the bracket equals $\sigma[\sinh(\alpha_p\tau-u)-\sinh\tau]$; since $\sinh$ is strictly increasing and $\alpha_p\tau-u<\tau$ by $\alpha_p<1$ and $u>0$, it is always negative. At $\theta=0$ the sign is that of $\sinh(\alpha_p\tau)/\sinh\tau-\sigma C$. This hyperbolic ratio strictly decreases from $\alpha_p$ to zero. If $\alpha_p>\sigma C$, there is a unique positive root $\tau_*$ at which the ratio equals $\sigma C$. For $0<\tau<\tau_*$ each $\tau$ gives a unique coupling angle in the fourth quadrant, and $\mathcal H_\theta>0$ gives real analyticity. This branch starts at $\tau=0,\theta=0$ and exits at $\tau=\tau_*,\theta=0$, with exit radius $e^{-2\tau_*/p}$. If $\alpha_p\le\sigma C$, the open candidate sector in the fourth quadrant contains no branch. If $p\ge2$, $\tau>0$ and $-\gamma<\theta<0$, then every term on the left of \eqref{eq:coupling-hyp} is nonpositive and the last one is strictly negative, so the open sector in the fourth quadrant contains no coupling line. These conclusions also appear in the first paper~\cite[Section 3.3 and Appendix A.1]{PartI}.

\paragraph{Half-angle equation at fixed radius.}

One may also set $t:=\tan(\theta/2)$ at fixed $r$. From
$\cos\theta=(1-t^2)/(1+t^2)$ and $\sin\theta=2t/(1+t^2)$, multiplying \eqref{eq:coupling-hyp} by the positive quantity $1+t^2$ gives
\begin{align*}
&-\big[S\sinh(\alpha_p\tau)+\sigma SC\sinh\tau\big]t^2
+2\delta C\cosh(\alpha_p\tau)t\\
&\hspace{8em}+S\sinh(\alpha_p\tau)-\sigma SC\sinh\tau=0 .
\end{align*}
For $p<2$ the coefficient of the quadratic term is strictly negative and that of the linear term strictly positive. When the constant term is positive the product of the two roots is negative, so the smaller root is the only negative one. By $\mathcal H_\theta>0$ in the fourth quadrant and $\mathcal H(\tau,-\gamma)<0$, this root corresponds to $-\gamma<\theta<0$. If the constant term vanishes, the smaller root becomes $t=0$, which corresponds to the exit on the real axis. When the constant term is negative, all three terms are strictly negative for any $t<0$, so there is no negative root. In the first quadrant one has to select $0<t<\tan(\gamma/2)$, and the number of roots is decided by the radial folding criterion above.

\section[\texorpdfstring{$1<p,q<2$ or $2<p,q<\infty$}{1<p,q<2 or 2<p,q<infty}]{$1<p,q<2$ or $2<p,q<\infty$}\label{sec:sameside}

\paragraph{Plan of this section.}
For $1<p,q<2$ the ridge line and the relevant coupling branch lie in the fourth quadrant. The first paper~\cite{PartI} provides the angular residual $\mathcal F$ along the coupling line and proves that $\mathcal F'=\mathcal K\mathcal N>0$ at each of its zeros. Each crossing therefore runs from negative to positive, so there is at most one intersection.

For $2<p,q<\infty$ we first map the stationary point to the dual problem with exponents $q',p'<2$ and then show that the resulting point lies in the fourth-quadrant candidate sector. This reduces the claim to uniqueness in the previous case. If $q\le2\le p$, no parameters outside the Weissler region exist; the proof is at the end of Section~\ref{sec:duality}.

\subsection{Proof of Theorem~\ref{thm:same} for \texorpdfstring{$1<p,q<2$}{1<p,q<2}}

By Section~\ref{sec:prelim} the ridge line is real analytic in $-\gamma<\theta<0$. If the sector contains a coupling line, that line is a single real analytic curve running from $z=1$ to the exit point on the real axis. The global maximum outside the Weissler region exceeds $1$ and is attained in the open candidate sector, so the two lines meet at least once. We now quote the intersection sign criterion of the first paper to prove uniqueness.

Parametrize the coupling line by $\tau=-\frac p2\log r$ and put
\[
w_{\rm c}:=-\frac CS\tan\theta,\qquad
t_{\rm c}:=\tanh\tau,\qquad u_{\rm c}:=\tanh(\alpha_p\tau),\qquad
\zeta_{\rm c}:=t_{\rm c}\frac{\delta u_{\rm c}-w_{\rm c}}{u_{\rm c}-\delta w_{\rm c}}.
\]
By~\cite[Section 3.5--3.8]{PartI}, the angular residual along the coupling line may be taken as
\[
\mathcal F(\tau):=\alpha_q\arctanh\zeta_{\rm c}-\arctanh w_{\rm c}.
\]
At a putative intersection point let $\eta_{\rm c}:=\arctanh\zeta_{\rm c}$. Then
$w_{\rm c}=\tanh(\alpha_q\eta_{\rm c})$, $0<w_{\rm c}<u_{\rm c}<1$ and $0<\zeta_{\rm c}<t_{\rm c}<1$. The derivative factorization of the first paper reads
\begin{equation}\label{eq:same-factor}
\mathcal F'(\tau)=\frac{\alpha_q\zeta_{\rm c}'}{1-\zeta_{\rm c}^2}
-\frac{w_{\rm c}'}{1-w_{\rm c}^2}=\mathcal K\mathcal N,
\end{equation}
where
\[
\mathcal K:=
\frac{(1-t_{\rm c}^2)(1-u_{\rm c}^2)^2(1-w_{\rm c}^2)^2}
{(u_{\rm c}^2-w_{\rm c}^2)
[t_{\rm c}w_{\rm c}(1-u_{\rm c}^2)(1-\zeta_{\rm c}^2)
+u_{\rm c}\zeta_{\rm c}(1-t_{\rm c}^2)(1-w_{\rm c}^2)]}>0.
\]
Writing the second factor with the three intersection quantities
\[
\begin{gathered}
P_{\rm c}:=\tanh\tau\,\sinh(\alpha_p\tau)\cosh(\alpha_p\tau),\quad
Q_{\rm c}:=\tanh\eta_{\rm c}\,\sinh(\alpha_q\eta_{\rm c})\cosh(\alpha_q\eta_{\rm c}),\\
D_{\rm c}:=\cosh^2(\alpha_p\tau)-\cosh^2(\alpha_q\eta_{\rm c})
\end{gathered}
\]
gives
\[
\begin{aligned}
\mathcal N:={}&D_{\rm c}^2+(\alpha_q-\alpha_p)(P_{\rm c}+Q_{\rm c})D_{\rm c}\\
&+\alpha_p\alpha_q(\cosh^2\tau-\cosh^2\eta_{\rm c})
\left(\frac{Q_{\rm c}^2}{\sinh^2\eta_{\rm c}}-\frac{P_{\rm c}^2}{\sinh^2\tau}\right)>0.
\end{aligned}
\]
The positivity of $\mathcal N$ is the content of the intersection sign criterion of~\cite[Section 3.8]{PartI} and holds at every interior intersection when $1<p,q<2$.

Every zero is therefore a strict upward crossing. If there were two zeros, the first would be positive just to its right and the second negative just to its left; between these two places there would be a zero first reached from the positive side, contradicting \eqref{eq:same-factor}. So the intersection is unique and is exactly the global maximizer.

\subsection{Duality and the case \texorpdfstring{$2<p,q<\infty$}{2<p,q<infty}}\label{sec:duality}

\paragraph{Correspondence of stationary points.}
Pairing with respect to the two-point probability measure gives $T_\omega^*=T_{\overline\omega}$, hence $N_{p,q}(\omega)=N_{q',p'}(\overline\omega)$. Let $J_a(f):=(\abs{f_j}^{a-2}f_j)_j$, with zero components mapped to zero.

\begin{lemma}[Dual correspondence of stationary points]\label{lem:dual-stationary}
Let $1<p,q<\infty$, $\omega\ne0$, and let $f$ be a projective stationary point of the quotient $f\mapsto\norm{T_\omega f}_q/\norm f_p$. Set $g:=J_q(T_\omega f)$. Then $g$ is a projective stationary point of the quotient $g\mapsto\norm{T_{\overline\omega}g}_{p'}/\norm g_{q'}$, and $f\mapsto g$ is a two-sided correspondence.
\end{lemma}
\begin{proof}
The projective stationarity equation is equivalent to
\begin{equation}\label{eq:euler}
T_\omega^* g=\frac{\norm{T_\omega f}_q^q}{\norm f_p^p}J_p(f).
\end{equation}
It comes from taking real derivatives of $\log(\norm{T_\omega f}_q/\norm f_p)$, and homogeneity allows all real variation directions. Denote the positive proportionality factor on the right by $\kappa$. The key identity is $J_{a'}(J_a(x))=x$, which is equivalent to $(a-1)(a'-2)+a-2=0$. Applying $J_{p'}$ to both sides of \eqref{eq:euler} and using $J_{p'}(\kappa y)=\kappa^{p'-1}J_{p'}(y)$ gives $J_{p'}(T_\omega^*g)=\kappa^{p'-1}f$; on the other hand $J_{q'}(g)=J_{q'}(J_q(T_\omega f))=T_\omega f$. Together these two relations are precisely the stationarity equation for $g$ in the dual problem, whose operator is $T_{\overline\omega}=T_\omega^*$ and whose source and target exponents are $q'$ and $p'$. Repeating the same procedure in the reverse direction recovers $f$, so the correspondence is a bijection. Moreover
\[
\norm{g}_{q'}=\norm{T_\omega f}_q^{q-1},\qquad
\norm{T_\omega^*g}_{p'}=\kappa\,\norm f_p^{p-1},
\]
so the two norm ratios agree and the correspondence preserves the value of the objective.
\end{proof}

\begin{lemma}[The dual input falls in the open sector]\label{lem:dual-sector}
Let $2<p,q<\infty$, let $\omega$ satisfy \eqref{eq:domain}, and let $z=re^{i\theta}$ be a stationary point of the quotient $\mathscr R$ with $0<r<1$ and $0<\theta<\gamma$. Write $w:=(c+z)/(1+cz)$. Then
\[
0<\arg w<\gamma .
\]
Consequently the input ratio $z'=\overline{\abs w^{\,q-2}w}$ of the dual problem satisfies $-\gamma<\arg z'<0$, that is, $z'$ lies in the fourth-quadrant candidate sector of $N_{q',p'}(\omega)$.
\end{lemma}
\begin{proof}
\textit{Step 1: $w$ lies in the right half-disc.} From
\[
\begin{aligned}
U-V&=(1-\rho^2)(1-r^2)+4\rho rS\sin\theta>0,\\
\operatorname{Re}\bigl[(c+z)\overline{(1+cz)}\bigr]
&=\operatorname{Re}c\,(1+r^2)+(1+\rho^2)\operatorname{Re}z>0
\end{aligned}
\]
we get $\abs w<1$ and $\operatorname{Re}w>0$. The latter implies $\arg w\in(-\pi/2,\pi/2)$. Also $\arg z'=-\arg w$, since multiplication by the positive number $\abs w^{\,q-2}$ does not change the argument.

\textit{Step 2: the argument is fixed by the dual angular equation.}
By Lemma~\ref{lem:dual-stationary} and complex conjugation, $z'$ is a projective stationary point for the parameter $\omega$ and the exponents $(q',p')$. Let $r':=|z'|\in(0,1)$ and $\theta':=\arg z'\in(-\pi/2,\pi/2)$, and take $U,V$ at $(r',\theta')$. Both are positive in this right half-disc, because $|z'|<1$ and $z'=-c$ is not in the right half-plane. By \eqref{eq:ridge}, the angular derivative has the sign of
\[
\sin(\gamma-\theta')-(V/U)^{p'/2-1}\sin(\gamma+\theta').
\]
For $0\le\theta'<\gamma$,
\[
U-V=(1-\rho^2)(1-r'^2)+4\rho r'\sin\gamma\sin\theta'>0 .
\]
Since $p'<2$ we have $(V/U)^{p'/2-1}>1$. Also $\sin(\gamma+\theta')\ge\sin(\gamma-\theta')>0$, so the angular derivative is strictly negative, in particular at $\theta'=0$. For $\gamma\le\theta'<\pi/2$ one sine is nonpositive and the other positive, so the derivative is again strictly negative. For $-\pi/2<\theta'\le-\gamma$ one sine is positive and the other nonpositive, so the derivative is strictly positive. The dual stationary point must therefore satisfy $-\gamma<\theta'<0$, whence $0<\arg w<\gamma$.
\end{proof}

\paragraph{Completion of the proof.}
Let $\omega$ lie outside the region. By Lemma~\ref{lem:argmax} the quotient $\mathscr R$ attains its global maximum in the open candidate sector, so that sector contains at least one stationary point. If two distinct stationary points existed, Lemma~\ref{lem:dual-stationary} and Lemma~\ref{lem:dual-sector} would produce two distinct stationary points of the quotient with parameter $\omega$ and exponents $(q',p')$ in its fourth-quadrant candidate sector; but $q',p'<2$, and Section~\ref{sec:sameside} shows that on the same side of $2$ at low exponents the candidate sector contains exactly one stationary point. This is a contradiction. Hence the quotient $\mathscr R$ has exactly one stationary point in the first-quadrant candidate sector, and it is the global maximizer. The two problems have the same norm, and the constant direction gives $1$ for both, so the outside condition corresponds across the duality. If $q\le2\le p$, monotonicity of the probability norms and the $L^2$ contraction of $T_\omega$ give directly
\[
\norm{T_\omega f}_q\le\norm{T_\omega f}_2\le\norm f_2\le\norm f_p .
\]
For $s<1$ the middle inequality is an equality only when all Fourier coefficients other than the constant one vanish, so there is no outside case here.

\section[\texorpdfstring{$1<p<2=q$ or $p=2<q<\infty$}{1<p<2=q or p=2<q<infty}]{$1<p<2=q$ or $p=2<q<\infty$}\label{sec:edges}

Let $1<p<2=q$. By the Parseval identity on two points,
\[
\RR(re^{i\theta})=
2^{1/p-1}\frac{\sqrt{(1+s^2)(1+r^2)+2(1-s^2)r\cos\theta}}{(1+r^p)^{1/p}}.
\]
For fixed $r>0$, since $1-s^2>0$, the unique maximal angle is $\theta=0$. From $2\operatorname{Re}c/(1+|c|^2)=(1-s^2)/(1+s^2)$ we get $\sigma C=(1-s^2)/(1+s^2)$. Maximizing the value on the positive real axis is therefore the same as maximizing
\[
\frac{1+r^2+2\sigma Cr}{(1+r^p)^{2/p}}.
\]
Differentiating and cancelling the positive factor $2(1+r^p)^{-2/p-1}$ leaves
\[
(r+\sigma C)(1+r^p)-(1+r^2+2\sigma Cr)r^{p-1}
=(1-r^p)\left(\sigma C-\frac{r^{p-1}-r}{1-r^p}\right).
\]
If $\xi=-\log r>0$, then
\[
\frac{r^{p-1}-r}{1-r^p}=
\frac{\sinh((2-p)\xi/2)}{\sinh(p\xi/2)}.
\]
By the same comparison of $t\coth t$ as before, this ratio strictly decreases with $\xi$ from $(2-p)/p$ to $0$, and therefore strictly increases with $r$ from $0$ to $(2-p)/p$. Hence outside the region the unique root $r_p(s)\in(0,1)$ is given by
\begin{equation}\label{eq:edge-root}
\frac{r^{p-1}-r}{1-r^p}=\sigma C
\end{equation}
and the profile on the positive real axis first increases and then decreases at this point. The root exists exactly when $\sigma C<(2-p)/p$, and cross multiplication turns this into $s^2>p-1$, which is precisely the condition $\omega\notin A_{p,2}$. Consequently $N_{p,2}(\omega)=\RR(r_p(s))$. Setting $\theta=0$ in \eqref{eq:coupling-hyp} gives the same root in $0<r<1$ as \eqref{eq:edge-root}, so the radius of this maximizer is also the exit radius on the real axis of the fourth-quadrant coupling branch. To count the stationary points of the whole projective space, note that the angular derivative vanishes only on the real axis. On the negative real axis the sign of the radial derivative is that of $r-r^{p-1}-\sigma C(1-r^p)$. This expression is negative for $0<r<1$ and positive for $r>1$. Together with the reciprocal symmetry and the fact that $\RR$ increases away from the origin along the positive real axis, the projective space has exactly $1$, $-1$ and one pair of reciprocal non-Fourier stationary points. On the positive real axis the point $1$ is a minimum in the radial direction and a maximum in the angular one. On the negative real axis the point $-1$ is simultaneously the minimal angle at fixed radius and the unique minimum along that axis, hence is the global minimizer.

If $p=2<q$, let $t:=r_{q'}(s)$, obtained from \eqref{eq:edge-root} with $q'$ in place of $p$. Duality gives $N_{2,q}(\omega)=N_{q',2}(\overline\omega)$, and the outside condition is $s^2>q'-1=1/(q-1)$. The dual extremal input is $(1,t)^\top$. The stationary point correspondence of the previous section then gives, when the original input exponent equals $2$, an original input proportional to $T_{\overline\omega}(1,t)^\top$, because $J_2$ is the identity map. Hence the unique maximizer in the disc is
\begin{equation}\label{eq:edge-recover}
z=\frac{\overline c+t}{1+t\overline c}.
\end{equation}
The squared modulus of the denominator minus that of the numerator is $(1-t^2)(1-\rho^2)>0$, so $\abs z<1$. After multiplying by the conjugate of the denominator, the imaginary part of the numerator is $(1-t^2)\rho\sin\gamma>0$ and the real part is $(1+t^2)\rho\cos\gamma+t(1+\rho^2)>0$. Dividing the two gives $\tan(\arg z)<\tan\gamma$, hence $0<\arg z<\gamma$. This also shows that the unique non-Fourier stationary point in the disc on this edge falls in the corresponding candidate sector.

\section[\texorpdfstring{$1<p<2<q<\infty$}{1<p<2<q<infty}]{$1<p<2<q<\infty$}\label{sec:count}

\paragraph{Plan of this section.}
For $1<p<2<q$ the candidate region is $0<r<1$, $0<\theta<\gamma$. Each circle carries a unique ridge angle. The coupling line is single-valued when the angle is used as a parameter, but may develop two branches when the radius is used. By Remark~\ref{rem:coordinate}, the ridge line and the coupling line in polar coordinates correspond pointwise to $\Phi_y=0$ and $\Phi_x=0$. We use the equality cases of H\"older's inequality to introduce two weights. For fixed weights, maximizing over the input yields the largest singular value of a weighted matrix. Taking the supremum over both weights gives $N_{p,q}(\omega)$. We then locate the intersections of the two curves in weight coordinates.

We first prove the classification for the dual exponent $p=q'$. In that case $P=Q$ and the objective is symmetric under exchanging the two weights. Comparing the layers at fixed $H$, and using the diagonal point coordinate $d$ together with the junction coordinate $u_c$, then decides whether there is one peak or two outside the region. In the general case we first use $\Phi_y=0$ to eliminate $y$ and then write $\Phi_x=0$ as $\mathcal P(x)=P$. Lemma~\ref{lem:elimination} gives the properties of the elimination curve and its two limits, Lemmas~\ref{lem:comparison} and~\ref{lem:turn} give the possible directions of folding away from the diagonal, and Section~\ref{sec:diag-analytic} handles the diagonal point. Finally the type of a stationary point is read off from the direction in which it crosses a horizontal line, and the two peak heights are compared by exchanging the two weights.

\subsection{Two diagonal weights and the largest singular value}\label{sec:weights}

Throughout this section $1<p<2<q$. Put
\begin{equation}\label{eq:weight-PQ}
P:=\frac{p}{2-p},\qquad Q:=\frac{q}{q-2},\qquad
A:=\frac{1+s^2}{2s},\qquad B:=\cos\phi .
\end{equation}
In particular $A>1$ and $0<B<1$. For $k>1$ define
\[
D_k(t):=\cosh(kt)^{-1/(2k)}
\begin{pmatrix}e^{t/2}&0\\0&e^{-t/2}\end{pmatrix}.
\]
Its diagonal vector has $L^{2k}$ norm $1$ with respect to the probability measure. As $t$ runs over the real axis the diagonal ratio runs over all positive numbers, so this gives all strictly positive two-point weights of $L^{2k}$ norm $1$. Weights with a zero component are obtained by a limit. Since $1/p=1/(2P)+1/2$, H\"older's inequality gives
\[
\norm f_p\le\norm{D_P(x)^{-1} f}_2.
\]
When both components of $f$ are nonzero, equality requires the weight to be proportional to $\abs f^{\,1-p/2}$, so the minimum equals $\norm f_p$. Likewise, $1/2=1/(2Q)+1/q$ gives $\norm{D_Q(y)g}_2\le\norm g_q$, with equality when the weight is proportional to $\abs g^{\,q/2-1}$, so the maximum equals $\norm g_q$. With zero components these identities hold in the sense of infima and suprema. When the data have no zero component, H\"older equality gives
\begin{equation}\label{eq:optimal-weights}
x=\frac{2-p}{2}\log\frac{\abs{f_1}}{\abs{f_2}},\qquad
y=\frac{q-2}{2}\log\frac{\abs{g_1}}{\abs{g_2}},
\qquad g=T_\omega f.
\end{equation}
Taking the supremum over the input and over both weights at the same time, and then writing $f=D_P(x)a$, therefore yields
\[
N_{p,q}(\omega)^2=\sup_{x,y\in\mathbb R}
\norm{D_Q(y)T_\omega D_P(x)}_{2\to2}^{\,2}.
\]
All suprema are taken successively over the same joint quotient, so the order may be exchanged.

Dropping the two $\cosh$ scalars, let $a=(1+\omega)/2$ and $b=(1-\omega)/2$. The four entries of the weighted matrix are, in order, $a e^{(x+y)/2},\, b e^{(y-x)/2},\, b e^{(x-y)/2},\, a e^{-(x+y)/2}$. The trace of their Gram matrix is
\[
2|a|^2\cosh(x+y)+2|b|^2\cosh(x-y)
=(1+s^2)\cosh x\cosh y+2sB\sinh x\sinh y=2sH,
\]
where
\begin{equation}\label{eq:weight-objective}
\begin{gathered}
H:=A\cosh x\cosh y+B\sinh x\sinh y,\qquad \Gamma:=\arcosh H,\\
\Phi_{P,Q}:=\Gamma-\frac{\log\cosh(Px)}{P}-\frac{\log\cosh(Qy)}{Q}.
\end{gathered}
\end{equation}
Both diagonal exponential weights have determinant $1$, so the Gram determinant equals $\abs{a^2-b^2}^2=s^2$. Also $H\ge A\cosh(|x|-|y|)\ge A>1$, so the two eigenvalues differ and the larger one is $s(H+\sqrt{H^2-1})=s e^\Gamma$. Restoring the normalization factors gives
\begin{equation}\label{eq:supphi}
N_{p,q}(\omega)^2=s\sup_{x,y\ge0} e^{\Phi_{P,Q}(x,y)}.
\end{equation}
We may take $x,y\ge0$: for fixed absolute values the same sign of the two weights maximizes $H$, and changing both signs at once leaves $H$ unchanged.

\begin{lemma}[Correspondence of stationary points and of local maximality]\label{lem:stationary}
Stationary points in the open first-quadrant candidate sector correspond one-to-one to stationary points of $\Phi_{P,Q}$ with $x,y>0$. At corresponding points
\begin{equation}\label{eq:recover}
r=\exp[-2x/(2-p)],\qquad
\tan\theta=\tan\gamma\tanh y,\qquad
\RR(re^{i\theta})^2=s\, e^{\Phi_{P,Q}(x,y)}.
\end{equation}
Local maximality, including its strictness, also corresponds. Inside the candidate sector a saddle point detected in the weight plane corresponds to a saddle point of the original quotient.
\end{lemma}

\begin{proof}
Consider the joint quotient in the three variables $x,y,z$:
\[
\frac{\norm{D_Q(y)T_\omega(1,z)^\top}_2^2}{\norm{D_P(x)^{-1}(1,z)^\top}_2^2}.
\]
With the input fixed, H\"older equality gives the unique optimal weights \eqref{eq:optimal-weights}. Write $g=T_\omega f$. In the candidate sector $0<|z|<1$ and $U>V>0$, so both weight variables are strictly positive. One also checks uniqueness of the stationary point of the joint quotient in the two weights. The derivative in $x$ of the logarithm of the denominator equals
$\tanh(Px)+\tanh(x-\log(|f_1|/|f_2|))$, which strictly increases. The derivative in $y$ of the logarithm of the numerator equals
$\tanh(y+\log(|g_1|/|g_2|))-\tanh(Qy)$, whose sign is decided by
$\log(|g_1|/|g_2|)-(Q-1)y$. The zeros of the two derivatives are exactly \eqref{eq:optimal-weights}. Since the partial derivatives of the joint quotient in the two weights vanish at the optimum, the terms that the chain rule produces through the two weights all disappear when one inserts the optimal weights and differentiates in $z$.

At fixed weights, use $a=D_P(x)^{-1}f$ to rewrite the joint quotient as the Rayleigh quotient of a Gram matrix. Its stationarity equation is $Ga=\lambda a$, because differentiating on the sphere in an arbitrary real tangential direction gives $Ga-\lambda a=0$. Ignoring a common positive scalar, its off-diagonal entry is
\[
G_{12}=\frac{1-s^2}{2}\cosh y-i\,\operatorname{Im}(\omega)\sinh y .
\]
Its real part is strictly positive and $G_{21}=\overline{G_{12}}$. Substituting $\lambda=G_{22}$ into the characteristic polynomial gives $\det(G-G_{22}I)=-|G_{12}|^2<0$, so the two eigenvalues lie on opposite sides of $G_{22}$. From the eigenvector equation in the second row,
\[
a_2/a_1=\frac{G_{21}}{\lambda-G_{22}} .
\]
A nonzero numerator forces both coordinates to be nonzero, and the right-hand side has positive real part for the larger eigenvalue and negative real part for the smaller one. The weights are positive and do not change the phase. Hence only the leading singular direction can occur in the candidate sector, and its phase satisfies
\[
\tan\theta=\frac{2\operatorname{Im}\omega}{1-s^2}\tanh y=\tan\gamma\tanh y .
\]
When we differentiate in the weights, the leading singular input is stationary with respect to the direction variable at fixed weight. Applying the chain rule once more, what remains after eliminating the input is the partial derivative of the joint quotient in the weights. The leading singular value is simple and admits a locally smooth choice given by the explicit eigenvector ratio, so this computation is legitimate.

Thus a stationary point of the original quotient produces, after the two eliminations, a positive stationary point of $\Phi_{P,Q}$. Conversely, recovering the leading singular direction from a stationary point with positive weights makes both partial derivatives of the joint quotient in the weights vanish. The first weight is given by the denominator formula above as $|z|=\exp[-2x/(2-p)]$, and the phase formula gives $0<\theta<\gamma$. Both output components are then nonzero and the second weight is given by the numerator formula above, so the original quotient is stationary as well. Substituting into the joint quotient gives the value in \eqref{eq:recover}. As a further check of the angular formula, the optimal output weight $2y=(q-2)\log(|g_1|/|g_2|)$ together with $|g_1|/|g_2|=\sqrt{U/V}$ gives
\[
e^{2y}=(U/V)^{(q-2)/2}=\frac{\sin(\gamma+\theta)}{\sin(\gamma-\theta)}
\]
at a stationary point. Expanding the two sines and using $(e^{2y}-1)/(e^{2y}+1)=\tanh y$ recovers the same phase formula.

Finally, each elimination attains its unique maximum at fixed values of the other group of variables, and the optimal point depends continuously on the retained variables. If the reduced objective has a strict local maximum, the original joint quotient is also strictly maximal near the corresponding point, because for every choice of the retained variables the joint quotient does not exceed its value at the optimum of that group. The converse follows by taking the continuous optimal values of the eliminated variables. Applying this argument to each elimination gives the equivalence of the two objectives. At a saddle point of the weight objective the original objective takes larger values in every neighborhood, since otherwise it would be a local maximum; moving the angle away from a nondegenerate angular maximum on the same circle makes it strictly smaller. Hence it is a saddle point too.
\end{proof}

The stationarity equations in the weight plane are therefore $\Phi_x=\Gamma_x-\tanh(Px)=0$ and $\Phi_y=\Gamma_y-\tanh(Qy)=0$.

\begin{remark}[Ridge line, coupling line and the two weight stationarity equations]\label{rem:coordinate}
In the corresponding region of candidate maximizers the two curves of Section~\ref{sec:prelim} can be transformed into $\Phi_y=0$ and $\Phi_x=0$, respectively. This uses the leading singular input coordinates at fixed weights. The first formula of \eqref{eq:recover} holds on the coupling line, while the radius of a general weight point is given by the invertible change of variables below.

For $r>0$ and $|\theta|<\gamma$, and in order to write the three coefficients of the weighted numerator, put
\[
\begin{gathered}
s_\pm(\theta):=\sin(\gamma\pm\theta),\qquad d_\gamma:=\rho\sin(2\gamma),\\
a_\theta:=s_+(\theta)+\rho^2s_-(\theta),\qquad
b_\theta:=\rho^2s_+(\theta)+s_-(\theta).
\end{gathered}
\]
All three coefficients are positive. Define
\begin{equation}\label{eq:polar-weights}
e^{2y}:=\frac{s_+(\theta)}{s_-(\theta)},\qquad
e^{2x}:=\frac{b_\theta+d_\gamma/r}{a_\theta+d_\gamma r}.
\end{equation}
The first formula is equivalent to $\tan\theta=\tan\gamma\,\tanh y$, and as a map from $(-\gamma,\gamma)$ to $\mathbb R$ it is strictly increasing and onto. With $\theta$ fixed, the numerator of the second formula strictly decreases and its denominator strictly increases, so the whole ratio drops strictly from $+\infty$ to $0$. Hence \eqref{eq:polar-weights} is an invertible change of variables from this polar region to the real weight plane, and it depends only on the matrix parameters.

To see where the change of variables comes from, fix $x,y$. After dropping positive factors independent of $r,\theta$, finding the leading singular input is the same as maximizing
\[
\frac{e^yU+e^{-y}V}{e^{-x}+e^xr^2}.
\]
In the numerator the coefficients of $\cos\theta$ and $\sin\theta$ are proportional to
$\cos\gamma(e^y+e^{-y})$ and $\sin\gamma(e^y-e^{-y})$, and the first of these is positive. Hence the maximal angle at fixed $r$ satisfies the first formula of \eqref{eq:polar-weights}. After inserting this angle the numerator is proportional to $a_\theta+2d_\gamma r+b_\theta r^2$. With the angle frozen, the radial derivative has the sign of
$b_\theta r+d_\gamma-e^{2x}r(a_\theta+d_\gamma r)$. By the strict monotonicity in \eqref{eq:polar-weights}, the unique radial root is a maximum, and it gives the leading singular input. This root is exactly the second formula.

When differentiating in the weights at this optimal input, the chain rule terms through $r,\theta$ all vanish. Taking the logarithm of the joint quotient and adding the two normalizations gives
\begin{equation}\label{eq:separate-weight}
\begin{aligned}
\Phi_y&=\tanh\left(y+\tfrac12\log(U/V)\right)-\tanh(Qy),\\
\Phi_x&=\tanh(-\log r-x)-\tanh(Px).
\end{aligned}
\end{equation}
Here $\Gamma_y$ and $\Gamma_x$ are the first $\tanh$ in the two lines. The input radius of a general weight point is therefore
\[
r=\exp[-x-\arctanh\Gamma_x(x,y)].
\]
Only on $\Phi_x=0$ can we substitute $\arctanh\Gamma_x=Px$ and so obtain \eqref{eq:recover}.

Since $\tanh$ strictly increases, the first line of \eqref{eq:separate-weight} vanishes exactly when
$e^{2y}=(U/V)^{(q-2)/2}$. Inserting the first formula of \eqref{eq:polar-weights} gives the angular stationarity equation \eqref{eq:ridge}. The second line vanishes exactly when $e^{2x}=r^{p-2}$. Inserting the second formula of \eqref{eq:polar-weights} gives
$b_\theta r+d_\gamma=r^{p-1}(a_\theta+d_\gamma r)$, which is precisely what one obtains from \eqref{eq:coupling-polar} after cross multiplication and cancellation of the common term $b_\theta r^{p+1}$. Both correspondences hold along the whole curve. In particular, solving $\Phi_y=0$ first and then intersecting it with $\Phi_x=0$ amounts, in polar coordinates, to following the ridge line and intersecting it with the coupling line.

Formula \eqref{eq:polar-weights} and the exponential equations of the two curves are valid for all $1<p,q<\infty$. When $q=2$ the ridge line is $y=0$, and when $p=2$ the coupling line is $x=0$. If $p,q\ne2$, one assigns to $P,Q$ the signs dictated by the algebraic formulas in \eqref{eq:weight-PQ}, and the correspondence of zero sets for \eqref{eq:separate-weight} persists. The argument of this section, which interprets the supremum of the weight objective as the norm, is used for $1<p<2<q$. In this range an intersection in the open candidate sector satisfies $r<1$ and $\theta>0$, so both weights at an intersection are strictly positive.
\end{remark}

\subsection{Layered extrema and peak counts at the dual exponent}\label{sec:dual-proof}

Let $p=q'$, so that $P=Q>1$. To state the criteria of this section, define for a general parameter $k>1$
\[
R_k(t):=\frac{\sinh(kt)}{\sinh t}\quad(t>0),\qquad R_k(0):=k .
\]
Also write $c_+:=(A+B)/2$ and $c_-:=(A-B)/2$.
Since $t\coth t$ strictly increases, $R_P'(t)>0$, and $R_P$ increases from $P$ to infinity. At the dual exponent $p-1=(P-1)/(P+1)$. Substituting $q=p/(p-1)$ into \eqref{eq:weissler} gives
\[
(2-p)\abs{p-1+\omega^2}\le p(p-1-s^2).
\]
The right-hand side must be nonnegative. Since $\abs{p-1+\omega^2}^2=(p-1-s^2)^2+4(p-1)s^2B^2$, squaring both sides and cancelling $4(p-1)>0$ gives
\[
(2-p)^2s^2B^2\le(p-1-s^2)^2.
\]
Under the sign condition already imposed on the right-hand side, taking square roots yields exactly $(2-p)sB\le p-1-s^2$. Rewriting with $\sqrt{A^2-1}=(1-s^2)/(2s)$ and $P=p/(2-p)$ gives $A+B\le P\sqrt{A^2-1}$.

\begin{theorem}[Peak count and norm at the dual exponent]\label{thm:dual}
Let $p=q'$ and suppose $\omega$ lies outside the Weissler region. Let $d$ and $u_c$ be the unique positive roots of the two equations
\begin{equation}\label{eq:dd}
R_P(d)\sqrt{A^2-1+(1-B^2)\tanh^2 d}=A+B,\qquad
R_P(u_c)=P\,\frac{A+B}{A-B}.
\end{equation}
If $2d\le u_c$, the open candidate sector contains exactly one stationary point, given by the weights $(d,d)$, and it attains the global maximum. If $2d>u_c$, there are exactly three stationary points; the weights $(d,d)$ give a saddle point, and the other two optimal weights are uniquely determined up to interchange by
\begin{equation}\label{eq:dualsystem}
\frac{R_P(u)}{R_P(v)}=\frac{c_+}{c_-},\qquad v>0,\qquad
H^2-m^2 W^2=1,
\end{equation}
where
\[
H=c_+\cosh u+c_-\cosh v,\qquad
W=\cosh(Pu)+\cosh(Pv),\qquad
m=\frac{c_+}{R_P(u)}=\frac{c_-}{R_P(v)}.
\]
Setting $(x,y)=((u+v)/2,(u-v)/2)$, or exchanging the two coordinates, recovers the two maximal points, which have equal height. The squared norm is $s e^{\Phi_{P,P}}$ at the corresponding point. Inside the Weissler region and on its boundary the candidate sector contains no stationary point, the norm equals $1$, and the only maximal direction is $z=1$.
\end{theorem}

\begin{proof}
On the diagonal $x=y=t$ one has $H=A+(A+B)\sinh^2t$ and $H_x=H_y=c_+\sinh(2t)$. The derivative of the diagonal profile is $2(\Gamma_x-\tanh(Pt))$. To compare the two quantities in it, use the strictly increasing function $a\mapsto a/\sqrt{1-a^2}$. A direct computation gives
\[
H^2-1-H_x^2
=A^2-1+(A^2-B^2)\sinh^2t
=\cosh^2t\,[A^2-1+(1-B^2)\tanh^2t].
\]
Hence the order relation between $\Gamma_x$ and $\tanh(Pt)$ is the same as that between
\[
\frac{(A+B)\sinh t}{\sqrt{A^2-1+(1-B^2)\tanh^2t}}
\quad\text{and}\quad \sinh(Pt).
\]
Dividing by $\sinh t>0$ gives the first equation of \eqref{eq:dd} and the corresponding sign of the derivative. The left-hand side strictly increases, from $P\sqrt{A^2-1}$ to infinity. Outside the region there is a unique positive root $d$, and the diagonal profile first rises and then falls. Inside the region and on its boundary the diagonal profile strictly decreases. Since $c_+>c_->0$, the second equation also has a unique positive root.

Put $u:=x+y$ and $v:=|x-y|$, so that the positive weight region is $u>v\ge0$. By the addition formula for hyperbolic functions,
\[
H=c_+\cosh u+c_-\cosh v,\qquad
\cosh(Px)\cosh(Py)=\frac{\cosh(Pu)+\cosh(Pv)}2.
\]
At fixed $H$, maximizing $\Phi_{P,P}$ is the same as minimizing $W=\cosh(Pu)+\cosh(Pv)$. Substitute $\cosh u=(H-c_-\cosh v)/c_+$. For $0<v<u$ the chain rule gives
\begin{equation}\label{eq:layer-deriv}
\frac{dW}{d(\cosh v)}
=P\left[R_P(v)-\frac{c_-}{c_+}R_P(u)\right].
\end{equation}
More precisely, differentiating $\cosh(Pt)$ with respect to $\cosh t$ gives $PR_P(t)$, while at fixed $H$ the derivative of $\cosh u$ with respect to $\cosh v$ equals $-c_-/c_+$. Expression \eqref{eq:layer-deriv} strictly increases with $\cosh v$, because $v$ increases while $u$ decreases and $R_P$ strictly increases. It extends continuously to $v=0$, where $R_P(0)=P$. At the other end $u=v$ the derivative is $P(1-c_-/c_+)R_P(v)>0$. Hence, if the derivative at $v=0$ is nonnegative, it is strictly positive for $v>0$ and the unique minimum is at $v=0$, including the case where the derivative vanishes. If the initial derivative is negative, there is exactly one interior zero, which gives the unique minimizer $(u,v)$ of $W$ at fixed $H$ and corresponds in $(x,y)$ to two exchanged points. The sign of the initial derivative changes at $u=u_c$, so the two branches are joined at $H=c_+\cosh u_c+c_-$.

On the interior branch $R_P(u)/R_P(v)=c_+/c_-$, so $u$ and $H$ strictly increase with $v$, and $u>v$. From
\[
dW=PR_P(u)\,d(\cosh u)+PR_P(v)\,d(\cosh v)
=\frac{P}{m}\,dH
\]
one gets $dW/dH=P/m$. On the diagonal branch $dW/dH=PR_P(u)/c_+$. As the interior branch approaches the junction, $u\to u_c$ and $v\to0$, and $m\to c_+/R_P(u_c)$, so the derivatives of the two branches match. This shows that the first derivative along the layered maximal profile also matches, since the profile is $\arcosh H-P^{-1}\log(W/2)$. On the interior branch that derivative is
\[
\frac1{\sqrt{H^2-1}}-\frac1{mW}.
\]
Both denominators are positive, so its sign equals the sign of $1-J$, where $J:=H^2-(mW)^2$. The condition $J=1$ on the interior branch is exactly the stationarity of the layered maximum with respect to $H$.

We now prove the monotonicity of the interior branch: when $J\le A^2$ one has $dJ/dv>0$, and as $v\to\infty$ the limit of $J$ exceeds $1$. Let
\[
\ell(t)=\frac{d}{dt}\log R_P(t)=P\coth(Pt)-\coth t .
\]
We already have $\ell(t)>0$. Also $\ell'(t)=\sinh^{-2}t-P^2\sinh^{-2}(Pt)>0$, because $R_P(t)>P$. It remains to show that $\ell(t)/\tanh(t/2)$ strictly decreases. After differentiating and multiplying by a positive factor this is equivalent to $\ell-\ell'\sinh t>0$, that is,
\[
P\coth(Pt)+\frac{P^2\sinh t}{\sinh^2(Pt)}>\coth(t/2).
\]
The right-hand side equals the value of the left-hand side at $P=1$, because $\coth t+1/\sinh t=\coth(t/2)$. Fix $t>0$, regard $P$ temporarily as a parameter increasing from $1$, and put $w=Pt$. The derivative of the left-hand side in $P$, multiplied by $\sinh^2 w$, equals
\[
\sinh w\cosh w-w+2P\sinh t(1-w\coth w).
\]
Since $1-w\coth w<0$ and $P\sinh t\le\sinh w$, this is at least
\[
\sinh w(\cosh w+2)-w(1+2\cosh w).
\]
The quantity vanishes at $w=0$ and has derivative $2\sinh w(\sinh w-w)>0$, hence is strictly positive. This proves the required strict decrease.

Along the interior branch write $\Delta=u-v>0$ and $Z=R_P(u)R_P(v)$, and put
\[
K_0=Z\sinh^2(\Delta/2)-\sinh^2(P\Delta/2).
\]
Since $A=m(R_P(u)+R_P(v))$, one has $A^2-J=4m^2(1-K_0)$. Explicitly, when $J/m^2$ is expanded the two square terms are
$R_P(u)^2\cosh^2u-\cosh^2(Pu)=R_P(u)^2-1$ and $R_P(v)^2-1$. The cross term divided by $2$ equals $Z\cosh\Delta-\cosh(P\Delta)$, because the subtracted part contains $\sinh(Pu)\sinh(Pv)=Z\sinh u\sinh v$. Subtracting this expansion of $J/m^2$ from $A^2/m^2=(R_P(u)+R_P(v))^2$, and using $\cosh t-1=2\sinh^2(t/2)$, yields the identity.

Differentiating the interior branch equation logarithmically gives $u'=\ell(v)/\ell(u)$ and $Z'=2\ell(v)Z$, where the primes denote differentiation in $v$. Also, from $\coth a-\coth b=-\sinh(a-b)/(\sinh a\sinh b)$,
\[
\sinh\Delta-\frac{P\sinh(P\Delta)}{Z}
=(\ell(u)-\ell(v))\sinh u\sinh v .
\]
Differentiating $K_0$ directly,
\[
K_0'=Z'\sinh^2(\Delta/2)
+\frac{Z}{2}(u'-1)
\left[\sinh\Delta-\frac{P\sinh(P\Delta)}{Z}\right].
\]
After substitution this becomes
\[
\frac{2\ell(u)}{Z}K_0'
=2\ell(u)\ell(v)(\cosh\Delta-1)
-(\ell(u)-\ell(v))^2\sinh u\sinh v.
\]
Since $\ell$ strictly increases while $\ell(t)/\tanh(t/2)$ strictly decreases,
\[
1<\frac{\ell(u)}{\ell(v)}<\frac{\tanh(u/2)}{\tanh(v/2)}.
\]
Apply the strict increase of $z+1/z-2$ for $z>1$ to this inequality. Using
\[
\frac{(\tanh(u/2)-\tanh(v/2))^2}{\tanh(u/2)\tanh(v/2)}
=\frac{4\sinh^2((u-v)/2)}{\sinh u\,\sinh v},
\]
one sees that $2\ell(u)K_0'/Z>0$ in the formula above. Moreover $m=c_-/R_P(v)$ strictly decreases. When $J\le A^2$ one has $K_0\le1$, so
\[
J'=-8mm'(1-K_0)+4m^2K_0'>0.
\]
At the far end, from $R_P(t)=e^{(P-1)t}(1-e^{-2Pt})/(1-e^{-2t})$ and the interior branch equation,
$\Delta\to\log(c_+/c_-)/(P-1)>0$, while $m\to0$ and $m^2Z=c_+c_-$. Therefore
\[
\lim_{v\to\infty}J
=A^2+4c_+c_-\sinh^2\left(\frac{\log(c_+/c_-)}{2(P-1)}\right)>1.
\]

At the junction, if $2d<u_c$ the derivative of the diagonal profile is negative, so $J>1$ at the starting point of the interior branch. If $2d=u_c$ then $J=1$ at the junction, and since $1<A^2$ the interior branch has $J>1$ immediately to the right of the starting point. It can never drop below $1$ afterwards, because every time $J=1$ one has $J'>0$. Neither case therefore produces an interior stationary point. If $2d>u_c$, the derivative of the diagonal profile at the junction is positive, so $J<1$ at the starting point of the interior branch. The limit at the far end exceeds $1$, so a root exists, and every root crosses $1$ from below, so there is only one. In this case $(d,d)$ is a strict maximum along the diagonal, but moving slightly off the diagonal within the level set of $H$ through that point decreases $W$ and increases $\Phi_{P,P}$, so it is a saddle point. A stationary point off the diagonal must first satisfy the first-order condition within a level set and then the first-order condition along the interior branch, so all of them have been listed. The global maximum outside the region lies in the candidate sector, and combining this with the symmetry and Lemma~\ref{lem:stationary} gives the conclusion about the maximizers.

Inside the region and on its boundary the diagonal profile decreases everywhere, so $J>1$ at the junction, and the same sign argument excludes interior stationary points. The contraction bound is given in the first paper, and the constant direction makes the norm equal to $1$. Since $\RR$ increases away from the origin along the positive real axis, and on the unit circle only $\pm1$ are two-dimensional stationary points, with $\RR(-1)=s<1$, the only maximal direction is $1$.
\end{proof}

\subsection{Elimination of the stationarity equations and the behavior at the two ends}

We keep the notation of \eqref{eq:weight-objective} and put
\[
\alpha:=\arctanh\Gamma_x,\qquad \beta:=\arctanh\Gamma_y,\qquad h:=\frac{H}{\sqrt{H^2-1}}.
\]
These symbols refer to the weight plane. For $x,y>0$ one has $H_x,H_y>0$. Since $H_{xx}=H_{yy}=H$, or directly by the difference of squares of hyperbolic functions,
\begin{align}
H^2-1-H_y^2&=(A^2-1)\cosh^2 x+(1-B^2)\sinh^2 x,\label{eq:hy-sq}\\
H^2-1-H_x^2&=(A^2-1)\cosh^2 y+(1-B^2)\sinh^2 y.\nonumber
\end{align}
Hence $0<\Gamma_x,\Gamma_y<1$. Differentiating $\Gamma_x=H_x/\sqrt{H^2-1}$ directly gives
\[
\Gamma_{xx}=h(1-\Gamma_x^2),\qquad \Gamma_{yy}=h(1-\Gamma_y^2).
\]
Dividing by the appropriate factor gives $\alpha_x=\beta_y=h>1$. Also $h=(1-H^{-2})^{-1/2}$, and $H$ strictly increases in $x$ and $y$, so $h$ strictly decreases in both variables.

\begin{analyticift}[{\cite[Section 2.3]{KrantzParks2002}}]
Let $F$ be real analytic in a neighborhood of $(x_0,y_0)\in\mathbb R^2$ with $F(x_0,y_0)=0$ and $\partial_y F(x_0,y_0)\ne0$. Then in some neighborhood of $(x_0,y_0)$ the solution set of $F(x,y)=0$ is the graph of a unique real analytic function $y=Y(x)$ with $Y(x_0)=y_0$.
\end{analyticift}

\begin{lemma}[Elimination from the two weight stationarity equations]\label{lem:elimination}
Fix $Q$ and $\omega$.
\begin{enumerate}
\item For every $x>0$ there is a unique $y=Y(x)>0$ with $\beta(x,y)=Qy$, and $Y$ is bounded and real analytic on the positive axis.
\item Define
\begin{equation}\label{eq:Pfunction}
\mathcal P(x):=\frac{\alpha(x,Y(x))}{x}.
\end{equation}
Stationary points in the open candidate sector correspond to positive roots of $\mathcal P(x)=P$. Along $Y$ the quotient $\tanh x/\tanh Y(x)$ strictly increases. If $Y(d)=d$ for some $d$, then such a point is unique, $x<Y(x)$ for $x<d$ and $x>Y(x)$ for $x>d$, and $\mathcal P(d)=Q$.
\item For all $x>0$, $\mathcal P(x)>1$, and $\lim_{x\to\infty}\mathcal P(x)=1$.
\item Outside the Weissler region one has $\mathcal P(x)>P$ for all sufficiently small $x>0$.
\end{enumerate}
\end{lemma}

\begin{proof}
\textit{(1) Existence, uniqueness, boundedness and real analyticity of the positive root.}
The equation $\beta=Qy$ is equivalent to $\Gamma_y=\tanh(Qy)$ and to the formula for $\sinh(Qy)$ obtained by inverting $\tanh$,
\[
\sinh(Qy)=H_y\Big/\sqrt{H^2-1-H_y^2}.
\]
By \eqref{eq:hy-sq}, after dividing by $\sinh y>0$, this becomes
\begin{equation}\label{eq:Y}
\frac{\sinh(Qy)}{\sinh y}
\sqrt{A^2-1+(1-B^2)\tanh^2x}
=A+B\frac{\tanh x}{\tanh y}.
\end{equation}
For fixed $x>0$ the left-hand side strictly increases and the right-hand side strictly decreases. As $y\downarrow0$ the right-hand side tends to infinity while the left-hand side stays finite; as $y\to\infty$ the left-hand side tends to infinity while the right-hand side stays finite, so there is exactly one root. Let $F$ be the left-hand side of \eqref{eq:Y} minus the right-hand side. The function $F$ is real analytic for $x,y>0$: $\sinh y$ and $\tanh y$ do not vanish, and the quantity under the square root is strictly positive. The left-hand side has positive derivative in $y$ and the right-hand side negative derivative, so $\partial_y F(x,Y(x))>0$. By the real analytic implicit function theorem the positive root is given near every $x>0$ by a real analytic function. Pointwise uniqueness makes these local functions agree on overlaps, so $Y$ is real analytic on all of $(0,\infty)$. Multiplying \eqref{eq:Y} by $\tanh y$, the right-hand side is at most $A+B$ while the left-hand side is at least $\sqrt{A^2-1}\sinh(Qy)/\cosh y$, which tends to infinity as $y\to\infty$. Hence $Y$ is bounded.

\noindent\textit{(2) Correspondence of stationary points and the diagonal intersection.}
Put $u:=\tanh x$ and $v:=\tanh Y(x)$. Take logarithms in \eqref{eq:Y} and differentiate in $\log u$. The left-hand side gives
\[
\frac{d\log R_Q(y)}{d\log v}\frac{d\log v}{d\log u}
+\frac{(1-B^2)u^2}{A^2-1+(1-B^2)u^2},
\quad y=\arctanh v.
\]
The right-hand side can be written as $\log(Av+Bu)-\log v$, whose derivative is
\[
\frac{Av}{Av+Bu}\frac{d\log v}{d\log u}+\frac{Bu}{Av+Bu}
-\frac{d\log v}{d\log u}
=\frac{Bu}{Av+Bu}\left(1-\frac{d\log v}{d\log u}\right).
\]
Rearranging gives
\begin{equation}\label{eq:log-slope}
\frac{d\log v}{d\log u}
=\frac{\dfrac{Bu}{Av+Bu}-\dfrac{(1-B^2)u^2}{A^2-1+(1-B^2)u^2}}
{\dfrac{d\log R_Q(y)}{d\log v}+\dfrac{Bu}{Av+Bu}}
<1.
\end{equation}
The first term of the denominator equals $\frac{v}{1-v^2}\frac{R_Q'(y)}{R_Q(y)}>0$, so the denominator exceeds $Bu/(Av+Bu)>0$ while the numerator is smaller than the same quantity. Hence $d\log(u/v)/d\log u>0$, and $\log u$ strictly increases with $x$. This gives uniqueness of the diagonal intersection and the relation between the two sides. Exchanging $x$ and $y$ gives the same property on the other stationarity curve. In particular, if two stationary points have $x_1<x_2$ then necessarily $y_1>y_2$. If two $y$ were equal, the uniqueness of the root after exchanging would force the two $x$ to be equal. If $y_1<y_2$, strict monotonicity on the two curves would require $u/v$ and $v/u$ to increase strictly at the same time, a contradiction.

Since $\Phi_y=\Gamma_y-\tanh(Qy)$ and $\Phi_x=\Gamma_x-\tanh(Px)$, the function $Y$ obtained from $\beta=Qy$ turns the second stationarity equation into $\mathcal P(x)=P$. At the diagonal intersection symmetry gives $\alpha=\beta$, hence $\mathcal P(d)=Q$. A positive diagonal point exists exactly when $Q\sqrt{A^2-1}<A+B$, which follows directly from the monotonicity of the diagonal equation \eqref{eq:dd}. Boundedness also ensures $\lim_{x\to\infty}\tanh x/\tanh Y(x)>1$, where the limit may be infinite, because for large enough $x$ the quotient already exceeds $1$ and strictly increases.

\noindent\textit{(3) The limit at the right end.}
From the symmetric version of \eqref{eq:hy-sq},
\[
\sinh\alpha(x,y)=
\frac{A\sinh x+B\cosh x\tanh y}
{\sqrt{A^2-1+(1-B^2)\tanh^2y}}.
\]
The denominator lies between two strictly positive constants. For large enough $x$ the numerator is bounded below by $A\sinh x$ and above by $A\sinh x+B\cosh x$, both of which are of the order of a fixed positive multiple of $e^x$. Using $\arsinh t=\log(t+\sqrt{1+t^2})$ gives $\alpha(x,Y(x))=x+O(1)$, hence $\mathcal P(\infty)=1$. On the other hand $\alpha_x=h>1$ and $\alpha(0,y)>0$, and integrating from $0$ at fixed $y$ gives $\alpha(x,y)>x$, hence $\mathcal P(x)>1$.

\noindent\textit{(4) The limit at the left end and the Weissler condition.}
Put $h_c:=A/\sqrt{A^2-1}$ and $k_c:=B/\sqrt{A^2-1}>0$. Since $Y$ is bounded, any sequence $x\to0$ has a subsequence along which $Y(x)$ converges. If the limit $y_0$ is positive, passing to the limit in \eqref{eq:Y} gives $R_Q(y_0)=h_c$. This equation has a unique positive root when $Q<h_c$ and no positive root when $Q\ge h_c$. If $Q<h_c$ and the limit were zero, the left-hand side would tend to $Q\sqrt{A^2-1}<A$ while the right-hand side is always at least $A$, a contradiction. Hence for $Q<h_c$ every subsequential limit equals this unique positive root, $Y(x)$ converges to it, and $\alpha(x,Y(x))$ tends to a positive number, so $\mathcal P(0^+)=+\infty$. If $Q\ge h_c$, every subsequential limit is $0$, so $Y(x)\to0$.

In the latter case $H$ equals $A$ at the origin, its first derivatives vanish, and $H_{xx}=H_{yy}=A$, $H_{xy}=B$. Expanding $H_x,H_y$, dividing by $\sqrt{H^2-1}$ and using $\arctanh t=t+O(t^3)$ therefore gives
\[
\alpha=h_c x+k_c y+O((x+y)^2),\qquad
\beta=k_c x+h_c y+O((x+y)^2).
\]
Substituting $y=Y(x)$ into $\beta=QY(x)$ and dividing both sides by $x+Y(x)$ gives
\[
(Q-h_c)\frac{Y(x)}{x+Y(x)}
=k_c\frac{x}{x+Y(x)}+O(x+Y(x)).
\]
If $Q=h_c$, then $x/(x+Y(x))\to0$ and hence $Y(x)/x\to\infty$. Dividing the expansion of $\alpha$ first by $Y(x)$, with remainder $O(Y(x))$, gives $\alpha(x,Y(x))/Y(x)\to k_c$, so $\mathcal P(0^+)=+\infty$. If $Q>h_c$, the display above gives $x/(x+Y(x))\to(Q-h_c)/(Q-h_c+k_c)$, that is, $Y(x)/x\to k_c/(Q-h_c)$. The expansion of $\alpha$ may now be divided by $x$, the remainder being $O(x)$, and one obtains
\[
\mathcal P(0^+)=h_c+\frac{k_c^2}{Q-h_c}.
\]
Outside the region this value is strictly larger than $P$. Indeed, $Q>h_c$ is equivalent to $(q-1)s^2<1$, so $D$ in \eqref{eq:local} is strictly negative. Being outside the region means that the quadratic form has a positive direction, so its determinant $ED-K^2<0$. Putting $P-h_c$ and $Q-h_c$ over a common denominator gives
\[
P-h_c=\frac{2(p-1-s^2)}{(2-p)(1-s^2)},\quad
Q-h_c=\frac{2(1-(q-1)s^2)}{(q-2)(1-s^2)},\quad
k_c=\frac{2a_0}{1-s^2}.
\]
From $E=s^2-(p-1)+(q-2)a_0^2$, $D=(q-1)s^2-1-(q-2)a_0^2$ and $b_0^2=s^2-a_0^2$,
\[
ED-K^2=(p-1-s^2)(1-(q-1)s^2)-(2-p)(q-2)a_0^2 .
\]
Therefore $(P-h_c)(Q-h_c)-k_c^2=4(ED-K^2)/[(2-p)(q-2)(1-s^2)^2]<0$, which is the required inequality at the left end.
\end{proof}

Outside the region, $\mathcal P(x)-P$ is positive at the left end and negative at the right end. The derivative criteria below control its turning points, and through them the number of intersections with a horizontal line. If one only wants to decide whether a diagonal point exists, the limit of the quotient $\tanh x/\tanh Y(x)$ at the left end can also be read off from the computation above. That limit is $0$ for $Q\le h_c$ and $(Q-h_c)/k_c$ for $Q>h_c$. The limit at the right end exceeds $1$. This is consistent with the fact that there is exactly one positive diagonal intersection when $Q<h_c+k_c=(A+B)/\sqrt{A^2-1}$. When the left limit equals $1$, strict monotonicity makes the quotient exceed $1$ everywhere on the positive axis.

\subsection{Two comparison inequalities away from the diagonal}

\begin{lemma}\label{lem:comparison}
For free variables $x>y>0$ one has $\alpha(x,y)>\beta(x,y)$ and
\begin{equation}\label{eq:comp}
\frac{\alpha(x,y)/x-h(x,y)}{\sinh^2\alpha(x,y)}
<
\frac{\beta(x,y)/y-h(x,y)}{\sinh^2\beta(x,y)}.
\end{equation}
\end{lemma}

\begin{proof}
$H_x-H_y=(A-B)\sinh(x-y)>0$. Dividing both by the same positive quantity $\sqrt{H^2-1}$ and applying $\arctanh$ gives $\alpha>\beta$.

For the second inequality fix $h>1$ at the current point and define
\[
g(t):=\arctanh(\tanh t/h),\qquad
b(t):=\arctanh\bigl((B/A)\tanh t\bigr).
\]
Since $H_x/H=(\tanh x+(B/A)\tanh y)/(1+(B/A)\tanh x\tanh y)$, the addition formula for $\tanh$ gives
\[
g(\alpha)=x+b(y),\qquad g(\beta)=y+b(x).
\]
Because $b$ is positive and strictly increasing, $b(y)/x<b(x)/y$. We show next that
\[
J_h(t):=\frac{t/g(t)-h}{\sinh^2 t}
\]
strictly decreases for $t>0$. A direct computation gives $g'(t)=h^{-1}/[1+(1-h^{-2})\sinh^2t]$. Write
\[
w_t(z):=\frac{1}{1+(1-h^{-2})\sinh^2(tz)}\quad (0\le z\le1).
\]
Integrating and changing variables gives $g(t)=t h^{-1}\int_0^1w_t(z)\,dz$. Putting $t/g(t)-h$ over a common denominator, the numerator is $h\int_0^1(1-w_t(z))\,dz$. Since $1-w_t=(1-h^{-2})\sinh^2(tz)w_t$,
\[
J_h(t)=h(1-h^{-2})
\frac{\int_0^1 [\sinh(tz)/\sinh t]^2 w_t(z)\,dz}{\int_0^1w_t(z)\,dz}.
\]
Write $f_t(z):=[\sinh(tz)/\sinh t]^2$ and let $d\mu_t(z):=w_t(z)\,dz/\int_0^1 w_t$. Thus $J_h(t)$ differs from the mean of $f_t$ with respect to the probability measure $\mu_t$ by a positive constant. For $0<z<1$ the function $f_t(z)$ strictly decreases in $t$, and for fixed $t$ it increases in $z$.

Take $t_2>t_1>0$. To compare the two probability measures, note first that
\[
\partial_t\log w_t(z)
=-2(1-h^{-2})\,z
\frac{\tanh(tz)}{1-h^{-2}\tanh^2(tz)} .
\]
For $0\le z\le1$ the quotient in the bracket increases in $z$, because its numerator increases while its denominator stays positive and decreases. Hence $\partial_t\log w_t(z)$ decreases in $z$. Integrating over $t\in[t_1,t_2]$ shows that $w_{t_2}(z)/w_{t_1}(z)$ also decreases in $z$.

With respect to the probability measure $\mu_{t_1}$ the function $f_{t_1}$ is increasing while $w_{t_2}/w_{t_1}$ is decreasing. The rearrangement inequality gives
\[
\int_0^1 f_{t_1}\frac{w_{t_2}}{w_{t_1}}\,d\mu_{t_1}
\le
\left(\int_0^1 f_{t_1}\,d\mu_{t_1}\right)
\left(\int_0^1 \frac{w_{t_2}}{w_{t_1}}\,d\mu_{t_1}\right).
\]
Concretely, subtracting the two functions of the two integration variables gives two quantities of opposite sign, so the double integral of their product is nonpositive. Expanding and using that $\mu_{t_1}$ has total mass $1$ gives exactly this inequality. Dividing by the second positive integral on the right turns the left-hand side into $\int f_{t_1}\,d\mu_{t_2}$. Using $f_{t_2}(z)<f_{t_1}(z)$ for $0<z<1$,
\[
\int f_{t_2}\,d\mu_{t_2}
<\int f_{t_1}\,d\mu_{t_2}
\le\int f_{t_1}\,d\mu_{t_1}.
\]
This proves directly that $J_h(t)$ strictly decreases.

From $g(\alpha)/x=1+b(y)/x$,
\[
\frac{\alpha/x-h}{\sinh^2\alpha}
=J_h(\alpha)+\frac{b(y)}x\left[J_h(\alpha)+\frac{h}{\sinh^2\alpha}\right].
\]
The quantity in brackets equals $\alpha/(g(\alpha)\sinh^2\alpha)>0$ and decreases strictly in $\alpha$. Exchanging $x$ and $y$ gives the corresponding formula on the right. Comparing term by term with $\alpha>\beta$ and $b(y)/x<b(x)/y$ yields \eqref{eq:comp}.
\end{proof}

\subsection{Direction of folding}

\begin{lemma}\label{lem:turn}
If $\mathcal P'(x_0)=0$ and $x_0\ne Y(x_0)$, then
\[
\sgn\mathcal P''(x_0)=\sgn(Y(x_0)-x_0).
\]
In other words, at $x_0<Y(x_0)$ the function can only turn from decreasing to increasing, and at $x_0>Y(x_0)$ only from increasing to decreasing.
\end{lemma}

\begin{proof}
Take $y_0:=Y(x_0)$ and $P_0:=\mathcal P(x_0)$, and fix the auxiliary exponents $(P_0,Q)$ throughout the computation. Then $(x_0,y_0)$ is a stationary point of $\Phi_{P_0,Q}$. All partial derivatives below, as well as $\alpha,\beta,h$, are evaluated at this point, so $P_0=\alpha/x_0$ and $Q=\beta/y_0$. From $\alpha_x=h$, $\alpha(0,y_0)>0$ and the strict decrease of $h$ in $x$,
\[
\alpha(x_0,y_0)=\alpha(0,y_0)+\int_0^{x_0}h(t,y_0)\,dt>x_0h(x_0,y_0).
\]
Hence $P_0-h>0$. Exchanging the variables gives $Q-h>0$ in the same way.

Define the following function, which is used to compare the two groups of terms in the third derivative along the tangent. For $\lambda>0$ let
\[
\mathscr K(\lambda):=
\lambda\,\frac{
3(h^2-1)-2(P_0-h)(Q-h)+2h(P_0-h)+(h^2-1)\lambda^2
}{
h^2-1+2h(P_0-h)+[3(h^2-1)-2(P_0-h)(Q-h)]\lambda^2
}.
\]
The proof has three steps, one for each relation.
\begin{enumerate}
\item $\sgn\mathcal P''(x_0)$ equals the sign of $\mathscr K\!\left(\sqrt{(P_0-h)/(Q-h)}\right)-\sinh\alpha/\sinh\beta$.
\item For every $\lambda>0$ the value $\mathscr K(\lambda)$ lies between $1$ and $\lambda$.
\item If $x_0>y_0$, then $\sinh\alpha/\sinh\beta$ strictly exceeds both $1$ and $\sqrt{(P_0-h)/(Q-h)}$. If $x_0<y_0$, it is strictly smaller than both.
\end{enumerate}
The three together give $\sgn\mathcal P''(x_0)=\sgn(Y(x_0)-x_0)$.

\medskip
\noindent\textit{Step 1. A sign formula from the third derivative along the tangent.}
Let $L_{P_0}(x):=\Phi_{P_0,Q}(x,Y(x))$. Along the curve $\Phi_y=0$, so
\[
L_{P_0}'(x)=\tanh(x\mathcal P(x))-\tanh(P_0x).
\]
At $x_0$ the two arguments coincide and both first derivatives equal $P_0$, so $L_{P_0}''(x_0)=0$. Differentiating once more, the second derivative terms of the outer $\tanh$ cancel. What remains is the second derivative of the inner function $x\mathcal P(x)$ at $x_0$, which equals $x_0\mathcal P''(x_0)$, so
\begin{equation}\label{eq:profile-third}
L_{P_0}'''(x_0)=x_0\sech^2\alpha\,\mathcal P''(x_0).
\end{equation}

Write $L_{P_0}''(x_0)=0$ and the derivative of $\Phi_y(x,Y(x))=0$ as
\[
\Phi_{xx}+\Phi_{xy}Y'(x_0)=0,\qquad
\Phi_{xy}+\Phi_{yy}Y'(x_0)=0 .
\]
Hence the tangent vector $\vect v:=(1,Y'(x_0))$ is annihilated by the Hessian, and $\Phi_{xx}=\Phi_{yy}Y'(x_0)^2$. The two diagonal entries are $-(P_0-h)\sech^2\alpha$ and $-(Q-h)\sech^2\beta$, so
\begin{equation}\label{eq:tangentrelation}
Y'(x_0)^2\frac{\sech^2\beta}{\sech^2\alpha}=\frac{P_0-h}{Q-h}.
\end{equation}
Writing the first stationarity equation backwards as $x=X(y)$, one has $X'(y_0)Y'(x_0)=1$ at that point, the slope relation of two tangent curves. If both slopes were positive, then after passing to the coordinates $\log\tanh x,\log\tanh y$ each slope is multiplied by a positive factor, and the two factors are reciprocal, so the product stays $1$. But \eqref{eq:log-slope} of Lemma~\ref{lem:elimination} and its exchanged version give that both positive slopes are strictly smaller than $1$, a contradiction. Hence $Y'(x_0)<0$. In applying \eqref{eq:log-slope} the first stationarity equation uses the auxiliary exponent $P_0$, which is the exponent appearing in the present Hessian.

Let $\vect c(x):=(x,Y(x))$, so that $\vect c'(x_0)=\vect v$. By the second-order chain rule,
\[
(\Phi_{P_0,Q}\circ\vect c)''=
D^2\Phi_{P_0,Q}[\vect c',\vect c']+D\Phi_{P_0,Q}[\vect c''] .
\]
Differentiate once more in $x$. In the first term, differentiating $D^2\Phi_{P_0,Q}(\vect c(x))$ produces an extra direction $\vect c'(x)$ from the variation of $\vect c(x)$, giving $D^3\Phi_{P_0,Q}[\vect c',\vect c',\vect c']$; differentiating the two input vectors produces one term $D^2\Phi_{P_0,Q}[\vect c',\vect c'']$ each. The derivative of the second term contributes one more such second-order term and one first-order term. Hence at $x_0$
\[
\begin{aligned}
L_{P_0}'''(x_0)
={}&D^3\Phi_{P_0,Q}[\vect v,\vect v,\vect v]\\
&+3D^2\Phi_{P_0,Q}[\vect v,(0,Y''(x_0))]
+D\Phi_{P_0,Q}[(0,Y'''(x_0))].
\end{aligned}
\]
Here $D^3\Phi_{P_0,Q}$ is a trilinear form and the first term expands as
\[
\Phi_{xxx}+3\Phi_{xxy}Y'(x_0)+3\Phi_{xyy}Y'(x_0)^2+\Phi_{yyy}Y'(x_0)^3.
\]
The other two terms vanish, since $D^2\Phi_{P_0,Q}\vect v=0$ and the gradient is zero. One may therefore compute along the fixed straight direction $\vect v$. Below $D_{\vect v}$ denotes differentiation in that direction. Since $H=\cosh\Gamma$,
\[
D_{\vect v}^3H
=\sinh\Gamma\bigl(D_{\vect v}^3\Gamma+(D_{\vect v}\Gamma)^3\bigr)
+3\cosh\Gamma\,D_{\vect v}\Gamma\,D_{\vect v}^2\Gamma .
\]
Substituting $D_{\vect v}\Gamma=\tanh\alpha+Y'(x_0)\tanh\beta$, and using $D_{\vect v}^2\Phi_{P_0,Q}=0$, gives
\[
D_{\vect v}^2\Gamma=P_0\sech^2\alpha+Y'(x_0)^2Q\sech^2\beta.
\]
Since $H_{xxx}=H_{xyy}=H_x$ and $H_{xxy}=H_{yyy}=H_y$,
\[
\frac{D_{\vect v}^3H}{\sinh\Gamma}
=(1+3Y'(x_0)^2)\tanh\alpha+(3Y'(x_0)+Y'(x_0)^3)\tanh\beta .
\]
Substituting these three expressions, and adding the third derivatives coming from the two negative logarithmic terms,
\[
2P_0^2\tanh\alpha\,\sech^2\alpha
+2Y'(x_0)^3Q^2\tanh\beta\,\sech^2\beta ,
\]
then expanding with $1-\tanh^2=\sech^2$ and grouping into four terms gives
\begin{align}
L_{P_0}'''(x_0)={}&
\tanh\alpha\,\sech^2\alpha(1-3hP_0+2P_0^2)\nonumber\\
&+3Y'(x_0)\tanh\beta\,\sech^2\alpha(1-hP_0)\nonumber\\
&+3Y'(x_0)^2\tanh\alpha\,\sech^2\beta(1-hQ)\nonumber\\
&+Y'(x_0)^3\tanh\beta\,\sech^2\beta(1-3hQ+2Q^2).
\label{eq:fourterms}
\end{align}
Dividing the first and third terms by $\tanh\alpha\,\sech^2\alpha>0$ and using \eqref{eq:tangentrelation} combines them into
\begin{align*}
&1-3hP_0+2P_0^2+3\frac{P_0-h}{Q-h}(1-hQ)\\
&=-\left[h^2-1+2h(P_0-h)
+[3(h^2-1)-2(P_0-h)(Q-h)]\frac{P_0-h}{Q-h}\right].
\end{align*}
Here we used $1-3hP_0+2P_0^2=1-h^2+h(P_0-h)+2(P_0-h)^2$. Dividing the second and fourth terms by the same positive quantity and factoring out $Y'(x_0)\tanh\beta/\tanh\alpha$, the remaining coefficients combine into
\begin{align*}
&3(1-hP_0)+\frac{P_0-h}{Q-h}(1-3hQ+2Q^2)\\
&=-\left[3(h^2-1)-2(P_0-h)(Q-h)+2h(P_0-h)+(h^2-1)\frac{P_0-h}{Q-h}\right].
\end{align*}
This time one expands $1-3hQ+2Q^2=1-h^2+h(Q-h)+2(Q-h)^2$. By \eqref{eq:tangentrelation} and $Y'(x_0)<0$,
\[
\frac{Y'(x_0)\tanh\beta}{\tanh\alpha}
=-\sqrt{\frac{P_0-h}{Q-h}}\frac{\sinh\beta}{\sinh\alpha}.
\]
The two groups therefore combine into exactly $\mathscr K\!\left(\sqrt{(P_0-h)/(Q-h)}\right)-\sinh\alpha/\sinh\beta$ times a positive factor. Precisely, that factor is the value of the first bracket divided by $\sinh\alpha/\sinh\beta$, and the bracket is positive by Step 2. Using \eqref{eq:profile-third} gives the first sign relation.

\medskip
\noindent\textit{Step 2. The value $\mathscr K(\lambda)$ lies between $1$ and $\lambda$.}
At the point of tangency the Hessian has zero determinant and $\Phi_{xy}=\Gamma_{xy}$, so
\[
(P_0-h)(Q-h)=\frac{\Gamma_{xy}^2}{\sech^2\alpha\,\sech^2\beta}.
\]
Expanding $H$ gives $HH_{xy}-H_xH_y=AB$, hence
$\Gamma_{xy}=(ABH-H_{xy})/(H^2-1)^{3/2}$. Moreover
\begin{equation*}
(H^2-1-H_x^2)(H^2-1-H_y^2)-(ABH-H_{xy})^2=(A^2-1)(1-B^2)(H^2-1).
\end{equation*}
To verify this identity, divide the left-hand side by $\cosh^2x\cosh^2y$, put $u:=\tanh x$ and $v:=\tanh y$, and use \eqref{eq:hy-sq} to get
\begin{align*}
&[A^2-1+(1-B^2)u^2][A^2-1+(1-B^2)v^2]\\
&\qquad-[B(A^2-1)-A(1-B^2)uv]^2\\
&=(A^2-1)(1-B^2)[A^2-1+u^2+v^2+2ABuv-(1-B^2)u^2v^2].
\end{align*}
The bracket in the last line equals $(A+Buv)^2-(1-u^2)(1-v^2)$, that is, $(H^2-1)/(\cosh^2x\cosh^2y)$. Therefore
\[
(P_0-h)(Q-h)=h^2-1-
\frac{(A^2-1)(1-B^2)}{(H^2-1-H_x^2)(H^2-1-H_y^2)}
<h^2-1.
\]
Consequently $3(h^2-1)-2(P_0-h)(Q-h)>h^2-1>0$, and the denominator of $\mathscr K$ is always strictly positive.

The numerator of $\mathscr K(\lambda)/\lambda-1$ is
\[
2[h^2-1-(P_0-h)(Q-h)](1-\lambda^2).
\]
The numerator of $\mathscr K(\lambda)-1$ is
\[
(\lambda-1)\big[(h^2-1)(\lambda-1)^2+2h(P_0-h)+2(P_0-h)(Q-h)\lambda\big].
\]
The latter is obtained by factoring out $\lambda-1$. The bracket is strictly positive, so $1\le\mathscr K(\lambda)\le\lambda$ for $\lambda\ge1$ and $\lambda\le\mathscr K(\lambda)\le1$ for $0<\lambda\le1$.

\medskip
\noindent\textit{Step 3. Comparison of $\sinh\alpha/\sinh\beta$ with the two endpoints.}
When $x_0>y_0$, Lemma~\ref{lem:comparison} gives $\alpha>\beta$ directly, hence $\sinh\alpha/\sinh\beta>1$. Substituting $P_0=\alpha/x_0$ and $Q=\beta/y_0$ into \eqref{eq:comp} gives
\[
\frac{P_0-h}{\sinh^2\alpha}<\frac{Q-h}{\sinh^2\beta}.
\]
Since $P_0-h>0$ and $Q-h>0$, rearranging and taking positive square roots gives $\sinh\alpha/\sinh\beta>\sqrt{(P_0-h)/(Q-h)}$. When $x_0<y_0$, applying the lemma to the exchanged coordinates reverses both comparisons. Combining the three relations gives the direction of folding.
\end{proof}

\subsection{Real analyticity at the diagonal point and absence of folding}\label{sec:diag-analytic}

By Lemma~\ref{lem:elimination}, the sets $x<Y(x)$ and $x>Y(x)$ occupy at most one interval each, with the former on the left. Lemma~\ref{lem:turn} shows that a horizontal tangency in the left interval is necessarily a nondegenerate local minimum, and a horizontal tangency in the right interval a nondegenerate local maximum. An interval cannot contain two such extrema. For instance, if two distinct strict local minima existed, one could take the maximum over the closed interval between them; near the higher endpoint the function still takes larger values, so the maximum would have to be attained inside, giving a contradictory local maximum. The right interval is treated in the same way. It remains to discuss \(x=d=Y(d)\).

If $\mathcal P'(d)\ne0$, there is no folding there. If $\mathcal P'(d)=0$, then $\mathcal P(d)=Q$ at this point, so we take the auxiliary objective $\Phi_{Q,Q}$. Differentiating along $y=Y(x)$,
\begin{equation}\label{eq:diag-identity}
\begin{aligned}
d\,\sech^2(Qd)\mathcal P'(d)
&=\left.\frac{d^2}{dx^2}\Phi_{Q,Q}(x,Y(x))\right|_{x=d}\\
&=\Phi_{xx}+\Phi_{xy}Y'(d)\\
&=\Phi_{yy}\bigl(1-Y'(d)^2\bigr)
=\bigl(1-Y'(d)\bigr)(\Phi_{xx}-\Phi_{xy}).
\end{aligned}
\end{equation}
The first line comes from differentiating $\tanh(x\mathcal P(x))-\tanh(Qx)$ and using $\mathcal P(d)=Q$. The rest uses $\Phi_{xx}=\Phi_{yy}$ and $\Phi_{xy}+\Phi_{yy}Y'(d)=0$. By \eqref{eq:log-slope}, the conversion factors of the two logarithmic coordinates coincide at $x=y=d$, so $Y'(d)<1$. Also $\Phi_{yy}<0$, so $\mathcal P'(d)=0$ forces $Y'(d)=-1$.

Exchange $x$ and $y$ in \eqref{eq:Y}. For each fixed $y>0$ the left-hand side strictly increases in $x$ and the right-hand side strictly decreases in $x$, and the two agree at $x=Y(y)$. Explicitly,
\[
\sgn\bigl(\Gamma_x(x,y)-\tanh(Qx)\bigr)=
\sgn\left[
 A+B\frac{\tanh y}{\tanh x}-
 R_Q(x)\sqrt{A^2-1+(1-B^2)\tanh^2y}
 \right].
\]
This is obtained by applying the strictly increasing function $a/\sqrt{1-a^2}$ to $\Gamma_x$ and $\tanh(Qx)$, and then multiplying by the positive factor
\[
\frac{\sqrt{A^2-1+(1-B^2)\tanh^2y}}{\sinh x}.
\]The bracket on the right is positive for $x<Y(y)$ and negative for $x>Y(y)$. Hence, taking $y=Y(x)$,
\begin{equation}\label{eq:composition-sign}
\sgn(\mathcal P(x)-Q)=\sgn\bigl(Y(Y(x))-x\bigr).
\end{equation}

The function $Y$ is real analytic on the positive axis and $Y(x)>0$, so the composition $x\mapsto Y(Y(x))$ is real analytic there as well. In detail, every $x_0>0$ has a real neighborhood on which the inner $Y$ equals a convergent power series; shrinking the neighborhood so that its image lies in the domain of analyticity of the outer $Y$ at $Y(x_0)$ gives local real analyticity of the composition. Let $F(x):=Y(Y(x))-x$. Then $F(d)=0$. If the zeros accumulated at $d$, all Taylor coefficients of $F$ at $d$ would vanish, since otherwise the first nonzero term would make $d$ an isolated zero. Thus $F$ would vanish identically near $d$, and by the real analytic identity theorem it would vanish on the whole connected positive axis. This contradicts the boundedness of $Y$ against the unboundedness of $x$. Hence $d$ is an isolated zero of $F$, and $F$ has a constant sign on each side.

Since $Y'(d)=-1$, the function $Y$ strictly decreases near $d$ and maps points on one side of $d$ to the other side. Take $x\ne d$ close enough to $d$ and let $x'=Y(x)$. Then $x$ and $Y(Y(x))$ both lie in the decreasing neighborhood. The quantity $F(x')=Y(Y(Y(x)))-Y(x)$ is obtained by applying the strictly decreasing $Y$ to the two terms of $F(x)=Y(Y(x))-x$ and subtracting, so the two have strictly opposite signs. They lie on opposite sides of $d$, so $F$ has opposite signs in the left and right neighborhoods of $d$. By \eqref{eq:composition-sign}, $\mathcal P-Q$ crosses zero at $d$. Moreover, the first nonzero Taylor term of $\mathcal P-Q$ at $d$ has odd degree. Hence $\mathcal P'$ has the same sign sufficiently close to $d$ on both sides. Thus, although $d$ may be a horizontal tangency, it does not change the direction of monotonicity of $\mathcal P$.

\subsection{Number of intersections, two-dimensional type and peak heights}\label{sec:peaks}

In summary, $\mathcal P$ has at most two turning points, and if there are two their order is necessarily a local minimum followed by a local maximum. The remaining pieces are strictly monotone, because $\mathcal P$ is real analytic on the positive axis, its derivative may vanish at $d$ but keeps its direction, and $\mathcal P>1$ together with $\mathcal P(+\infty)=1$ excludes a constant. Hence any horizontal line has at most three distinct intersections. Outside the region the opposite signs at the two ends give at least one. If there are three, the three pieces are decreasing, increasing and decreasing in this order, and the horizontal line lies strictly between the two turning heights. The three roots are therefore a downward crossing, an upward crossing and a downward crossing.

To read off the two-dimensional type from the direction of crossing, fix $x>0$ first. The same change of variables as in \eqref{eq:Y} gives
\[
\sgn\Phi_y(x,y)=\sgn\left[
 A+B\frac{\tanh x}{\tanh y}
 -R_Q(y)\sqrt{A^2-1+(1-B^2)\tanh^2x}
 \right].
\]
The bracket strictly decreases and vanishes at $y=Y(x)$, so at fixed $x$ the objective first strictly increases and then strictly decreases in $y$, and
\[
\Phi_{P,Q}(x,y)\le L_P(x):=\Phi_{P,Q}(x,Y(x)),
\]
with equality only at $y=Y(x)$. The derivative of the profile,
\begin{equation}\label{eq:profile}
L_P'(x)=\tanh(x\mathcal P(x))-\tanh(Px),
\end{equation}
has the sign of $\mathcal P(x)-P$. A downward crossing gives a strict maximum of the profile, and for nearby $x\ne x_0$ the upper bound above satisfies $L_P(x)<L_P(x_0)$. Changing $y$ at $x=x_0$ also strictly lowers the objective, so the point is a strict two-dimensional maximum. An upward crossing gives a strict minimum of the profile; moving along the profile away from the stationary point increases the objective, while changing $y$ at fixed $x$ decreases it, so the point is a saddle point. If there are exactly two distinct roots, one is a crossing and the other is a tangency of the horizontal line with some turning point. At a tangency $\mathcal P-P$ has the same sign on both sides, and $\mathcal P''\ne0$ at the turning point, so \eqref{eq:profile} and \eqref{eq:profile-third} show that the corresponding profile has a horizontal inflection point. A global maximum outside the region must be attained in the open candidate sector, so with one root that root is the global maximizer, and with two roots the downward crossing is the global maximizer.

If there are three roots, ordered as $x_1<x_2<x_3$, then $x_1$ is the outer peak and $x_3$ is the inner peak, because the radius in \eqref{eq:recover} strictly decreases in $x$. Let $x_-,x_+$ be the abscissas of the minimal and the maximal turning points. Then
\[
x_1<x_-<d<x_+<x_3 .
\]
Hence the outer peak satisfies $x_1<Y(x_1)$ and the inner peak satisfies $x_3>Y(x_3)$. Exchanging the two weights leaves $H$ unchanged, and
\begin{equation}\label{eq:swapheight}
\Phi_{P,Q}(x,y)-\Phi_{P,Q}(y,x)
=\int_y^x[\tanh(Qt)-\tanh(Pt)]\,dt .
\end{equation}
If $P<Q$, then for $x<y$ the exchange strictly increases the objective, so the outer peak cannot be the global maximizer. The middle root is a saddle point, so the only global peak is the inner one. If $P>Q$, the same argument excludes the inner peak. If $P=Q$, the middle root is $d$ and the two off-diagonal peaks are obtained from each other by the exchange, hence have equal height. Finally,
\[
P<Q\quad\Longleftrightarrow\quad
p(q-2)<q(2-p)
\quad\Longleftrightarrow\quad p<q'
\]
with the corresponding equality and reversed inequality. This gives the comparison of peak heights in Theorem~\ref{thm:main}.

\subsection{Deciding the folding interval by the exchange curvature}\label{sec:curvature}

Fix $Q$ and $\omega$. If $\mathcal P$ has two turning points, the positive diagonal point $d$ must lie between them, so the two turning values $P_-:=\mathcal P(x_-)$ and $P_+:=\mathcal P(x_+)$ satisfy $P_-<\mathcal P(d)=Q<P_+$. The horizontal line $Q$ meets the curve at $d$, and also on the right descending piece because $\mathcal P(+\infty)=1<Q$. By the symmetry of the coordinates under exchanging the two weights when $P=Q$, exchanging this off-diagonal intersection produces another one. Conversely, if the horizontal line $Q$ has three intersections, the curve must have two turning points. Before a positive diagonal point exists there certainly cannot be two turning points.

At the positive diagonal point $(d,d)$ the second derivative in the exchange direction is $2(\Phi_{xx}-\Phi_{xy})$, so by \eqref{eq:diag-identity} it has the sign of $\mathcal P'(d)$, because $d\sech^2(Qd)>0$ and $Y'(d)<1$. Using $H_x=H_y$ and $H-H_{xy}=A-B$,
\begin{equation}\label{eq:exchange}
\Phi_{xx}-\Phi_{xy}=
\frac{A-B}{\sqrt{H^2-1}}-Q\sech^2(Qd).
\end{equation}
From $\Gamma_x(d,d)=\tanh(Qd)$ and \eqref{eq:hy-sq},
\[
\sqrt{H^2-1}=\cosh d\,\cosh(Qd)\sqrt{A^2-1+(1-B^2)\tanh^2d}.
\]
Thus \eqref{eq:exchange} is positive exactly when
\[
(A-B)\cosh(Qd)>Q\cosh d\sqrt{A^2-1+(1-B^2)\tanh^2d}.
\]
By the diagonal point equation the square root equals $(A+B)\sinh d/\sinh(Qd)$. Multiplying by $\sinh(Qd)>0$ and using $\sinh(2t)=2\sinh t\cosh t$ gives
\begin{equation}\label{eq:threshold}
R_Q(2d)>Q\,\frac{A+B}{A-B}.
\end{equation}
Equality and the reversed inequality correspond in the same way.

If $\mathcal P'(d)>0$, then $\mathcal P(d)=Q$ while $\mathcal P>Q$ just to the right of $d$ and $\mathcal P<Q$ for large enough $x$, so there is another root on the right. Exchanging the two weights gives a third root. If $\mathcal P'(d)<0$, then $d$ would have to lie on an increasing piece whenever two turning points exist, so three roots are excluded directly. By Theorem~\ref{thm:dual}, the case $\mathcal P'(d)=0$ also belongs to the situation without two turning points. Below we give a perturbative argument by contradiction that does not use the layered proof for the dual line, in order to explain why equality cannot preserve three roots.

Suppose then that $\mathcal P'(d)=0$ while the horizontal line $Q$ has three roots. Since $Y'(d)=-1$, the slope of the curve $y=Y(x)$ at the diagonal point is negative. Both off-diagonal roots are downward crossings, and if one of them had zero derivative, Lemma~\ref{lem:turn} would force it to be a turning point, a contradiction. So both off-diagonal roots are simple.

Fix $A$ and $d$, increase $B$ slightly from its original value $B_0$, and define
\[
Q_B:=\beta_B(d,d)/d,\qquad P=Q_B.
\]
Here the subscript $B$ means that the quantity is computed with the modified $B$. In this way $(d,d)$ remains a stationary point of the corresponding dual problem. From
\[
R_{Q_B}(d)=\frac{A+B}{\sqrt{A^2-1+(1-B^2)\tanh^2 d}},
\]
the right-hand side strictly increases in $B$ while $R_{Q_B}(d)$ at fixed $d$ strictly increases in $Q_B$, so $Q_B$ increases. At the same time $H(d,d)=A\cosh^2d+B\sinh^2d$ increases, so $h_B(d,d)$ decreases. Hence $Q_B-h_B(d,d)>0$ strictly increases. From $\beta_B(x,Y_B(x))=Q_BY_B(x)$,
\[
Y_B'(d)=\frac{\partial_x\beta_B(d,d)}{Q_B-h_B(d,d)}.
\]
Here $\partial_x\beta_B=(ABH-H_{xy})/[\sqrt{H^2-1}(H^2-1-H_y^2)]$. At $x=y=d$,
\[
H_{xy}-ABH=A(1-B^2)\sinh^2d-B(A^2-1)\cosh^2d.
\]
Dividing numerator and denominator by $1-B^2>0$ gives
\[
-\partial_x\beta_B(d,d)=
\frac{ A\sinh^2d-\dfrac{B(A^2-1)}{1-B^2}\cosh^2d}
{\bigl[\dfrac{A^2-1}{1-B^2}\cosh^2d+\sinh^2d\bigr]
 \sqrt{(A\cosh^2d+B\sinh^2d)^2-1}} .
\]
At $B=B_0$ one has $Y_{B_0}'(d)=-1$, so the numerator of this fraction is positive. After a slight increase of $B$ the numerator is still positive but strictly smaller, because $B/(1-B^2)$ strictly increases. Both factors of the denominator strictly increase, so $-\partial_x\beta_B(d,d)$ strictly decreases. Together with the increase of $Q_B-h_B(d,d)$ this gives $-1<Y_B'(d)<0$. Now \eqref{eq:diag-identity} gives
\[
d\,\sech^2(Q_Bd)\mathcal P_B'(d)
=\Phi_{yy}(d,d)\bigl(1-Y_B'(d)^2\bigr)<0.
\]

Finally, $Q_B$ depends smoothly on $B$ by definition. In \eqref{eq:Y}, use the parameters $B,Q_B$ and subtract the right-hand side from the left. The resulting function has a nonzero partial derivative in $y$ at the solution. Thus near each original root, $Y_B(x)$ and $\mathcal P_B(x)-Q_B$ depend smoothly on $B$ and $x$. The two off-diagonal roots were simple, so by the implicit function theorem they both persist after the perturbation. Together with the point $d$, which is always retained, the perturbed problem still has three roots. But three roots force $d$ to lie on an increasing piece, contradicting $\mathcal P_B'(d)<0$. Hence equality cannot give three intersections.

\subsection{Classification of parameters and recovery of the norm}\label{sec:recover}

Given $p,q,\omega$ with $\omega\notin A_{p,q}$, the number of stationary points can be computed as follows.
\begin{enumerate}
\item Fix $q$ and $\omega$. If $Q\sqrt{A^2-1}\ge A+B$, there is no positive diagonal point. If $Q\sqrt{A^2-1}<A+B$, find the unique positive root $d$ of the first equation in \eqref{eq:dd}, with $R_Q$ in place of $R_P$, since the auxiliary exponent $P=Q$ is the one examined first.
\item If there is no positive diagonal point, or if $d$ exists but \eqref{eq:threshold} fails, there is at most one turning point. Outside the region the function $\mathcal P(x)-P$ has opposite signs at the two ends, so there is exactly one root. Here equality in \eqref{eq:threshold} also gives one root.
\item If \eqref{eq:threshold} holds, solve $\mathcal P'=0$ in $(0,d)$ and in $(d,\infty)$ to obtain the unique $x_-,x_+$. They can also be bracketed first by the three roots of the dual horizontal line $Q$: the left turning point lies between the left root and $d$, and the right turning point between $d$ and the right root. Put $P_\pm:=\mathcal P(x_\pm)$ and $p_\pm:=2P_\pm/(1+P_\pm)$. Since $p\mapsto p/(2-p)$ strictly increases, the classification is as follows.
\end{enumerate}
\begin{center}
\small
\begin{tabular}{llp{7.0cm}}
\toprule
$p$ with $\omega\notin A_{p,q}$ & stationary points & global maximizer\\
\midrule
$p<p_-$ or $p>p_+$ & 1 & the unique stationary point\\
$p=p_-$ & 2 & inner peak\\
$p_-<p<p_+$ & 3 & inner peak if $p<q'$; outer peak if $p>q'$; both if $p=q'$\\
$p=p_+$ & 2 & outer peak\\
\bottomrule
\end{tabular}
\end{center}
The first equality case $p=p_-$ is a tangency of the horizontal line at the left turning point, which merges the former outer peak with the middle saddle point. The other root lies on the right descending piece and corresponds to the inner peak. The case $p=p_+$ is the reverse. The derivative needed for the turning points can be written directly from \eqref{eq:Pfunction}:
\[
\mathcal P'(x)=\frac1x\left[
\frac{\Gamma_{xx}+\Gamma_{xy}Y'(x)}{1-\Gamma_x^2}-\mathcal P(x)\right],
\qquad Y'(x)=\frac{\beta_x}{Q-h},
\]
where the quantities on the right are evaluated at $(x,Y(x))$ and $Q-h>0$. Once it is known which root gives the global maximum, solve $\mathcal P(x)=P$ on the corresponding monotone piece, recover $z$ and $N_{p,q}(\omega)$ from \eqref{eq:recover}, and then obtain the induced norm of the original matrix from \eqref{eq:matrix-scaling}.

\subsection{A strict gap in the ordinary Hilbert factorization}\label{sec:gap}

Let $1<p<2<q$ and let $\omega$ lie outside the region. Take a set of optimal positive weights in \eqref{eq:supphi}. Replacing $\omega=se^{i\phi}$ by $s$ leaves $A$ unchanged while $B$ rises from $\cos\phi$ to $1$. At the same weights $H$ strictly increases, hence the largest singular value strictly increases. Therefore $N_{p,q}(\omega)<N_{p,q}(s)$. For any factorization $\mu\nu=\omega$ of the complex parameter, the two edge norms $N_{p,2}(\mu)$ and $N_{2,q}(\nu)$ depend only on the moduli of their arguments. For the first, the Parseval formula of Section~\ref{sec:edges} shows it: when the modulus exceeds $1$, the maximal angle at fixed radius is simply $\pi$. For the second it follows by duality. Hence
\[
N_{p,2}(\mu)N_{2,q}(\nu)
=N_{p,2}(|\mu|)N_{2,q}(|\nu|)
\ge N_{p,q}(s)>N_{p,q}(\omega),
\]
where the middle inequality comes from $T_{|\nu|}T_{|\mu|}=T_s$. The weighted singular value formula, by contrast, retains $\phi$ and recovers the exact norm together with the direction of the extremal input.

\clearpage
\subsection{A numerical plot of one folding interval}

Take $q=3$, $\rho=19/20$ and $\tan(\gamma/2)=9/10$. From $\cos\gamma=(1-\tan^2(\gamma/2))/(1+\tan^2(\gamma/2))$ and $\sin\gamma=2\tan(\gamma/2)/(1+\tan^2(\gamma/2))$, and using $\omega=(1-c)/(1+c)$, one obtains
\[
\omega=\frac{2353}{50727}+\frac{15200}{16909}i .
\]
For each given $x$, solve for the unique positive root $Y(x)$ of \eqref{eq:Y} by bisection, and then bracket the turning points on the two sides of $d$ with the derivative formula of Section~\ref{sec:recover}. The table keeps eight decimal places.
\begin{center}
\small
\begin{tabular}{lccc}
\toprule
 & $x$ & $\mathcal P(x)$ & $2\mathcal P(x)/(1+\mathcal P(x))$\\
\midrule
left turning point & 0.12479358 & 1.88508399 & 1.30677928\\
right turning point & 0.43345842 & 4.50757273 & 1.63686362\\
\bottomrule
\end{tabular}
\end{center}
At $p=q'=3/2$ the first coordinates of the three weights are approximately $0.02967283$, $0.30378373$ and $1.04452349$, corresponding to input radii approximately $0.88808188$, $0.29666998$ and $0.01532769$. By \eqref{eq:recover}, the value of $\RR$ at the two peaks is approximately $1.06877544$, and at the middle saddle point approximately $1.05452874$.

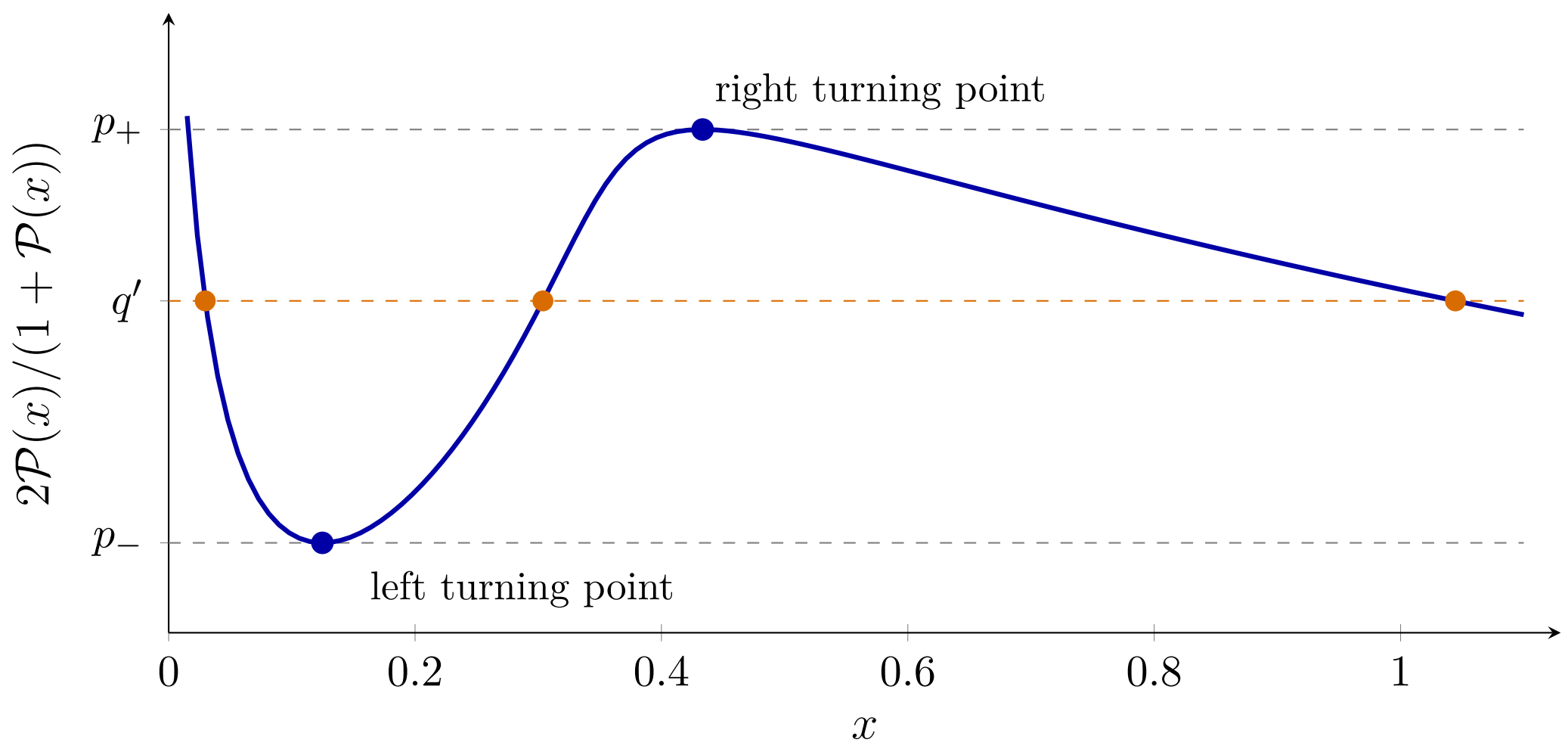
\begin{figure}[h!]
\centering
\begin{tikzpicture}
\begin{axis}[
width=.86\linewidth,height=7.0cm,
xlabel={$x$},ylabel={$2\mathcal P(x)/(1+\mathcal P(x))$},
xmin=0,xmax=1.13,ymin=1.235,ymax=1.73,
xtick={0,.2,.4,.6,.8,1.0},ytick={1.30677928,1.5,1.63686362},
yticklabels={$p_-$,$q'$,$p_+$},axis lines=left,
clip=false, tick label style={font=\small},
]
\addplot[gray,dashed,thin,forget plot] coordinates {(0,1.30677928) (1.1,1.30677928)};
\addplot[gray,dashed,thin,forget plot] coordinates {(0,1.63686362) (1.1,1.63686362)};
\addplot[orange!85!black,dashed,thin,forget plot] coordinates {(0,1.5) (1.1,1.5)};
\addplot[blue!65!black,very thick,no marks] coordinates {
(0.015000,1.647645) (0.023282,1.552383) (0.031565,1.486965) (0.039847,1.439791) (0.048130,1.404689) (0.056412,1.378061) (0.064695,1.357664) (0.072977,1.342017) (0.081260,1.330106) (0.089542,1.321204) (0.097824,1.314780) (0.106107,1.310435) (0.114389,1.307862) (0.122672,1.306822) (0.130954,1.307129) (0.139237,1.308631) (0.147519,1.311208) (0.155802,1.314763) (0.164084,1.319218) (0.172366,1.324509) (0.180649,1.330585) (0.188931,1.337405) (0.197214,1.344938) (0.205496,1.353160) (0.213779,1.362051) (0.222061,1.371601) (0.230344,1.381801) (0.238626,1.392647) (0.246908,1.404138) (0.255191,1.416275) (0.263473,1.429056) (0.271756,1.442479) (0.280038,1.456533) (0.288321,1.471193) (0.296603,1.486414) (0.304885,1.502111) (0.313168,1.518149) (0.321450,1.534316) (0.329733,1.550308) (0.338015,1.565727) (0.346298,1.580114) (0.354580,1.593028) (0.362863,1.604142) (0.371145,1.613318) (0.379427,1.620605) (0.387710,1.626185) (0.395992,1.630303) (0.404275,1.633211) (0.412557,1.635141) (0.420840,1.636284) (0.429122,1.636800) (0.437405,1.636815) (0.445687,1.636428) (0.453969,1.635717) (0.462252,1.634745) (0.470534,1.633561) (0.478817,1.632203) (0.487099,1.630704) (0.495382,1.629087) (0.503664,1.627375) (0.511947,1.625582) (0.520229,1.623724) (0.528511,1.621811) (0.536794,1.619854) (0.545076,1.617859) (0.553359,1.615834) (0.561641,1.613783) (0.569924,1.611713) (0.578206,1.609627) (0.586489,1.607527) (0.594771,1.605419) (0.603053,1.603303) (0.611336,1.601183) (0.619618,1.599059) (0.627901,1.596935) (0.636183,1.594810) (0.644466,1.592687) (0.652748,1.590567) (0.661031,1.588450) (0.669313,1.586338) (0.677595,1.584231) (0.685878,1.582130) (0.694160,1.580035) (0.702443,1.577948) (0.710725,1.575867) (0.719008,1.573795) (0.727290,1.571731) (0.735573,1.569675) (0.743855,1.567628) (0.752137,1.565590) (0.760420,1.563561) (0.768702,1.561541) (0.776985,1.559531) (0.785267,1.557531) (0.793550,1.555540) (0.801832,1.553559) (0.810115,1.551588) (0.818397,1.549628) (0.826679,1.547677) (0.834962,1.545736) (0.843244,1.543806) (0.851527,1.541886) (0.859809,1.539976) (0.868092,1.538076) (0.876374,1.536186) (0.884656,1.534307) (0.892939,1.532438) (0.901221,1.530579) (0.909504,1.528730) (0.917786,1.526892) (0.926069,1.525063) (0.934351,1.523245) (0.942634,1.521436) (0.950916,1.519638) (0.959198,1.517850) (0.967481,1.516072) (0.975763,1.514303) (0.984046,1.512545) (0.992328,1.510796) (1.000611,1.509057) (1.008893,1.507328) (1.017176,1.505609) (1.025458,1.503899) (1.033740,1.502199) (1.042023,1.500509) (1.050305,1.498827) (1.058588,1.497156) (1.066870,1.495494) (1.075153,1.493841) (1.083435,1.492197) (1.091718,1.490562) (1.100000,1.488937)
};
\addplot[only marks,mark=*,mark size=2.6pt,blue!65!black] coordinates {(0.12479358,1.30677928) (0.43345842,1.63686362)};
\addplot[only marks,mark=*,mark size=2.4pt,orange!85!black] coordinates {(0.02967283,1.50000000) (0.30378373,1.50000000) (1.04452349,1.50000000)};
\node[anchor=north west,font=\footnotesize] at (axis cs:.15,1.297) {left turning point};
\node[anchor=south west,font=\footnotesize] at (axis cs:.43,1.640) {right turning point};
\end{axis}
\end{tikzpicture}
\caption{The stationarity curve for $q=3$ and the matrix parameters chosen above. The orange horizontal line is the dual exponent $p=q'$, and the three orange dots are, from left to right, the outer peak, the saddle point and the inner peak. The two blue dots determine $p_-,p_+$.}
\label{fig:turning}
\end{figure}

\clearpage

\section[\texorpdfstring{norms at $(p,q)=(4/3,4)$ and $(3/2,3)$}{norms at two special exponents}]{The norms at $(p,q)=(4/3,4)$ and $(3/2,3)$}\label{app:closed}

The norms in this section are the two-point probability norms $N_{p,q}(\omega)$; the passage to counting measure is given in \eqref{eq:matrix-scaling}. The matrix parameters continue to satisfy \eqref{eq:domain} and to lie outside the region. The matrix parameters used here are $A=(1+s^2)/(2s)$ and $B=\cos(\arg\omega)$, where $s=|\omega|$. The closed formulas at these two exponents can also serve as a check on the general peak count.

\subsection{\texorpdfstring{$p=4/3,\ q=4$}{p=4/3, q=4}}

\begin{theorem}[\texorpdfstring{$p=4/3,\ q=4$}{p=4/3, q=4}]\label{thm:closed-four}
Let $J_4:=3A^2+2AB+3B^2$. Under the outside hypothesis of this section,
\begin{equation}\label{eq:norm4}
N_{4/3,4}(\omega)^2=
\begin{cases}
s\sqrt{\dfrac{2}{2-A^2-B^2}},&J_4<4,\\[5pt]
\dfrac{s(A+B)}{1-(A-B)^2/4},&J_4\ge4.
\end{cases}
\end{equation}
In the closed unit disc, there are exactly two global maximizers when $J_4<4$, obtained from each other by exchanging the weight coordinates. When $J_4\ge4$ there is exactly one global maximizer in the disc, corresponding to the positive diagonal weights. The two expressions agree at $J_4=4$.
\end{theorem}

\begin{proof}

Here $P=Q=2$. To obtain a compact closed form for the norm, put $u:=\arctan(\tanh x)$ and $v:=\arctan(\tanh y)$ in \eqref{eq:supphi}; for positive weights $0<u,v<\pi/4$. Write
\[
X_4:=\cos(u+v),\qquad Y_4:=\cos(u-v).
\]
Let $L_4$ denote the squared largest singular value at the given weights divided by $s$. The trace of the Gram matrix before $\cosh$ normalization, divided by $s$, equals $2H$, and its determinant divided by $s^2$ equals $1$. After normalization the trace is multiplied by $(\cosh(2x)\cosh(2y))^{-1/2}$ and the determinant by $(\cosh(2x)\cosh(2y))^{-1}$. Now
\[
\frac{\cosh x}{\sqrt{\cosh(2x)}}=\cos u,\quad
\frac{\sinh x}{\sqrt{\cosh(2x)}}=\sin u,\quad
\frac1{\cosh(2x)}=\cos(2u),
\]
and the same formulas hold for $y,v$. They follow directly from $\cosh(2x)=\cosh^2x+\sinh^2x$ and $\tan u=\tanh x$. Hence the normalized trace divided by $s$ equals
\[
2A\cos u\cos v+2B\sin u\sin v=(A-B)X_4+(A+B)Y_4.
\]
Also, from $\cos(2u)\cos(2v)=\frac12[\cos(2u+2v)+\cos(2u-2v)]=X_4^2+Y_4^2-1$, the determinant term simplifies as well. Thus $L_4$ is the larger root of
\[
L_4^2-L_4[(A-B)X_4+(A+B)Y_4]+X_4^2+Y_4^2-1=0 .
\]
Completing the square separately in the linear terms of $X_4$ and $Y_4$ gives
\[
\left(X_4-\frac{L_4(A-B)}2\right)^2+
\left(Y_4-\frac{L_4(A+B)}2\right)^2
=1-L_4^2\left(1-\frac{A^2+B^2}{2}\right).
\]

If $J_4<4$, then $A^2+B^2<2$. Since the right-hand side is nonnegative,
$L_4\le L_{4,\max}:=\sqrt{2/(2-A^2-B^2)}$. Take
\[
X_4=\frac{L_{4,\max}(A-B)}2,\qquad Y_4=\frac{L_{4,\max}(A+B)}2.
\]
Then $0<X_4<Y_4$, and $Y_4<1$ is equivalent to $J_4<4$. Also
$X_4^2+Y_4^2=(A^2+B^2)/(2-A^2-B^2)>1$, because $A>1$. Hence $a=\arccos X_4$ and $b=\arccos Y_4$ satisfy $0<b<a<\pi/2$ and $a+b<\pi/2$, the latter because $X_4^2+Y_4^2>1$ is equivalent to $\cos a>\sin b$. Therefore $u+v=a$ and $u-v=\pm b$ give exactly two pairs of admissible positive weights. Substituting back into the characteristic polynomial, the other root is $L_{4,\max}(A^2+B^2-1)<L_{4,\max}$, so the leading singular value is indeed attained. Both squares must vanish at the upper bound, so there are exactly two maximizers in the disc.

If $J_4\ge4$, choose any global maximizer in the open candidate sector. By Lemma~\ref{lem:stationary} it corresponds to positive weights. If $u\ne v$, the Jacobian of $(u,v)\mapsto(X_4,Y_4)$ does not vanish, because neither of the two sine factors $\sin(u+v)$, $\sin(u-v)$ is zero. Hence both partial derivatives of $L_4$ in $X_4,Y_4$ vanish. The leading singular value is a simple root, so implicit differentiation of the polynomial gives $2X_4=L_4(A-B)$ and $2Y_4=L_4(A+B)$ directly. Substituting back into the completed square gives $L_4=L_{4,\max}$, whose feasibility forces $J_4<4$, a contradiction. The optimal weights must therefore satisfy $u=v$, that is, $Y_4=1$. Differentiating along this line in the same way gives $2X_4=L_4(A-B)$, and substituting into the quadratic, then cancelling the positive factor $L_4$, gives $L_4[1-(A-B)^2/4]=A+B$. This determines $L_4$ and $X_4$ uniquely, hence also the positive diagonal weights. This yields \eqref{eq:norm4}. At $J_4=4$ one has $4=2(A^2+B^2)+(A+B)^2=(A-B)^2+2(A+B)^2$, so both expressions equal $2s/(A+B)$. The general dual criterion also simplifies directly here. Since $R_2(d)=2\cosh d$, put $t:=\sinh^2d$; squaring the first equation of \eqref{eq:dd} and eliminating $1+t$ gives $4[A^2-1+(A^2-B^2)t]=(A+B)^2$. Hence
\[
\cosh(2d)=1+2t
=1+\frac{(A+B)^2-4(A^2-1)}{2(A^2-B^2)}
=\frac{4-(A-B)^2}{2(A^2-B^2)}.
\]
Also $\cosh u_c=(A+B)/(A-B)$, and the numerator of the difference of the two expressions is $4-J_4$, so $2d>u_c$ corresponds exactly to $J_4<4$.

\end{proof}

\subsection{\texorpdfstring{$p=3/2,\ q=3$}{p=3/2, q=3}}

\begin{theorem}[\texorpdfstring{$p=3/2,\ q=3$}{p=3/2, q=3}]\label{thm:closed-three}
Let $t$ be the unique positive root of the cubic \eqref{eq:cubic} and $d=\operatorname{arsinh}\sqrt t$. If $R_3(2d)\le3(A+B)/(A-B)$, the squared norm is given by the single-peak formula \eqref{eq:norm3-single}. If $R_3(2d)>3(A+B)/(A-B)$, the squared norm is given by the two-peak formula \eqref{eq:norm3-double}, where $m$ is the larger root of \eqref{eq:quadratic}. The corresponding maximizers are recovered from \eqref{eq:recover}.
\end{theorem}

\begin{proof}

Here $P=Q=3$. By the triple angle formula,
\[
R_3(t)=3+4\sinh^2t=4\cosh^2t-1,\quad
\cosh(3t)=\cosh t\,[R_3(t)-2].
\]
Put $t:=\sinh^2d$. Substituting $R_3(d)=3+4t$ into \eqref{eq:dd}, where both sides are positive, squaring and multiplying by $1+t$ gives
\begin{equation}\label{eq:cubic}
(3+4t)^2[A^2-1+(A^2-B^2)t]=(A+B)^2(1+t).
\end{equation}
Write this as $(3+4t)^2=(A+B)^2(1+t)/[A^2-1+(A^2-B^2)t]$; the left-hand side strictly increases and the sign of the derivative of the right-hand side is decided by $B^2-1<0$, so the right-hand side strictly decreases. At $t=0$ the two sides equal $9$ and $(A+B)^2/(A^2-1)$ respectively, so outside the region there is exactly one positive root. This root is simple, since the derivative of the difference of the two sides is strictly positive. It can be found by bisection. For an explicit expression, the coefficients of the cubic \eqref{eq:cubic} from the cubic term down to the constant term are
\[
\begin{split}
a_3&=16(A^2-B^2),\\
a_2&=16(A^2-1)+24(A^2-B^2),\\
a_1&=24(A^2-1)+9(A^2-B^2)-(A+B)^2,\\
a_0&=9(A^2-1)-(A+B)^2 .
\end{split}
\]
Putting $t=z-a_2/(3a_3)$ gives $z^3+\mathfrak a z+\mathfrak b=0$, where
\[
\mathfrak a=\frac{3a_3a_1-a_2^2}{3a_3^2},\qquad
\mathfrak b=\frac{2a_2^3-9a_3a_2a_1+27a_3^2a_0}{27a_3^3},\qquad
\Delta=\frac{\mathfrak b^2}{4}+\frac{\mathfrak a^3}{27}.
\]
If $\Delta\ge0$, the unique positive root is given by the real cube roots
\[
t=\sqrt[3]{-\mathfrak b/2+\sqrt\Delta}
+\sqrt[3]{-\mathfrak b/2-\sqrt\Delta}-\frac{a_2}{3a_3}
\]
since the product of the two cube roots is $-\mathfrak a/3$, and expanding the cube of their sum verifies the original equation. If $\Delta=0$ and there is another multiple root, the formula picks the simple root, and the positive root has already been shown to be simple. If $\Delta<0$, then $\mathfrak a<0$; substituting $z=2\sqrt{-\mathfrak a/3}\cos\vartheta$ turns the equation into
$\cos(3\vartheta)=-\mathfrak b/[2(-\mathfrak a/3)^{3/2}]$. Of the trigonometric solutions the largest is obtained from $$\vartheta=\frac13\arccos(-\mathfrak b/[2(-\mathfrak a/3)^{3/2}]),$$
and this gives the unique positive root. Monotonicity thus determines the rule for selecting the root, after which radicals or trigonometric functions give the explicit solution.

When $2d\le u_c$, that is, $R_3(2d)\le3(A+B)/(A-B)$, the maximizer corresponds to the diagonal weights. From $\cosh^2(3d)=(1+t)(1+4t)^2$,
\begin{equation}\label{eq:norm3-single}
N_{3/2,3}(\omega)^2=
\frac{s\left(A+(A+B)t+\sqrt{(A+(A+B)t)^2-1}\right)}
{\big[(1+t)(1+4t)^2\big]^{1/3}}.
\end{equation}

When $2d>u_c$, we pass to the interior branch and use the parameter $m=c_+/R_3(u)=c_-/R_3(v)$. Since $v>0$, one has $0<m<c_-/3$. By the triple angle formula,
\[
\cosh^2u=\frac{c_++m}{4m},\quad
\cosh^2v=\frac{c_-+m}{4m},\quad
mW=(c_+-2m)\cosh u+(c_--2m)\cosh v.
\]
Hence $H-mW=2m(\cosh u+\cosh v)$ and
$H+mW=2[(c_+-m)\cosh u+(c_--m)\cosh v]$. Multiplying, the two square terms are $c_+^2-m^2$ and $c_-^2-m^2$. The sum of the two cross terms equals
$4m(A-2m)\cosh u\cosh v=(A-2m)\sqrt{(c_++m)(c_-+m)}$. The condition $J=1$ is therefore equivalent to
\begin{equation}\label{eq:quadratic-unsq}
1-c_+^2-c_-^2+2m^2=(A-2m)\sqrt{(c_++m)(c_-+m)} .
\end{equation}
The right-hand side is strictly positive. When squaring, the square of the right-hand side is
\[
(A-2m)^2(m^2+Am+c_+c_-)
=4m^4+(4c_+c_--3A^2)m^2+AB^2m+A^2c_+c_-,
\]
where the cubic terms cancel and $A^2-4c_+c_-=B^2$ was used. The square of the left-hand side is
$4m^4+4(1-c_+^2-c_-^2)m^2+(1-c_+^2-c_-^2)^2$. Subtracting, the quartic terms cancel and one obtains
\begin{equation}\label{eq:quadratic}
(4-B^2)m^2-AB^2m+\frac{(2-A^2-B^2)^2-A^2(A^2-B^2)}{4}=0 .
\end{equation}
The physical solution must lie in $0<m<(A-B)/6$ and satisfy \eqref{eq:quadratic-unsq}. It is exactly the larger root of the quadratic. Indeed, at the physical solution the difference of the two squares can be written as the sum of the two sides of \eqref{eq:quadratic-unsq} times the difference of the two sides, $1-J$. The first factor is strictly positive. The derivative of the second in $m$ is positive, because $dJ/dv>0$ at $J=1$ while $dm/dv<0$. Hence the quadratic has positive derivative at the physical solution, which identifies the larger root. Therefore
\[
m=\frac{AB^2+\sqrt{A^2B^4-(4-B^2)[(2-A^2-B^2)^2-A^2(A^2-B^2)]}}
{2(4-B^2)}.
\]
Finally, from $mW=\sqrt{H^2-1}$ the squared norm equals $s(H+mW)/(W/2)^{1/3}$. Substituting $c_\pm=(A\pm B)/2$ into $H+mW$ and $W/2$, and factoring $1/(2\sqrt{2m})$ out of the numerator and $1/(8m\sqrt{2m})$ out of the cube root in the denominator, the two factors cancel and leave $2^{-1/3}$. Thus
\begin{equation}\label{eq:norm3-double}
N_{3/2,3}(\omega)^2=
\frac{s}{\sqrt[3]{2}}\,
\frac{(A+B-2m)\sqrt{A+B+2m}+(A-B-2m)\sqrt{A-B+2m}}
{\big[(A+B-4m)\sqrt{A+B+2m}+(A-B-4m)\sqrt{A-B+2m}\big]^{1/3}}.
\end{equation}
Here $s=A-\sqrt{A^2-1}$. The parameters of the two branches are determined by \eqref{eq:cubic} and \eqref{eq:quadratic} respectively. The weights of the two peaks are recovered from
\[
u=\arcosh\sqrt{\frac{c_++m}{4m}},\qquad v=\arcosh\sqrt{\frac{c_-+m}{4m}},
\]
after which \eqref{eq:recover} gives the original input ratio.
\end{proof}

\clearpage
\appendix
\section{Plots of the one- and two-branch coupling lines}\label{app:coupling-plots}

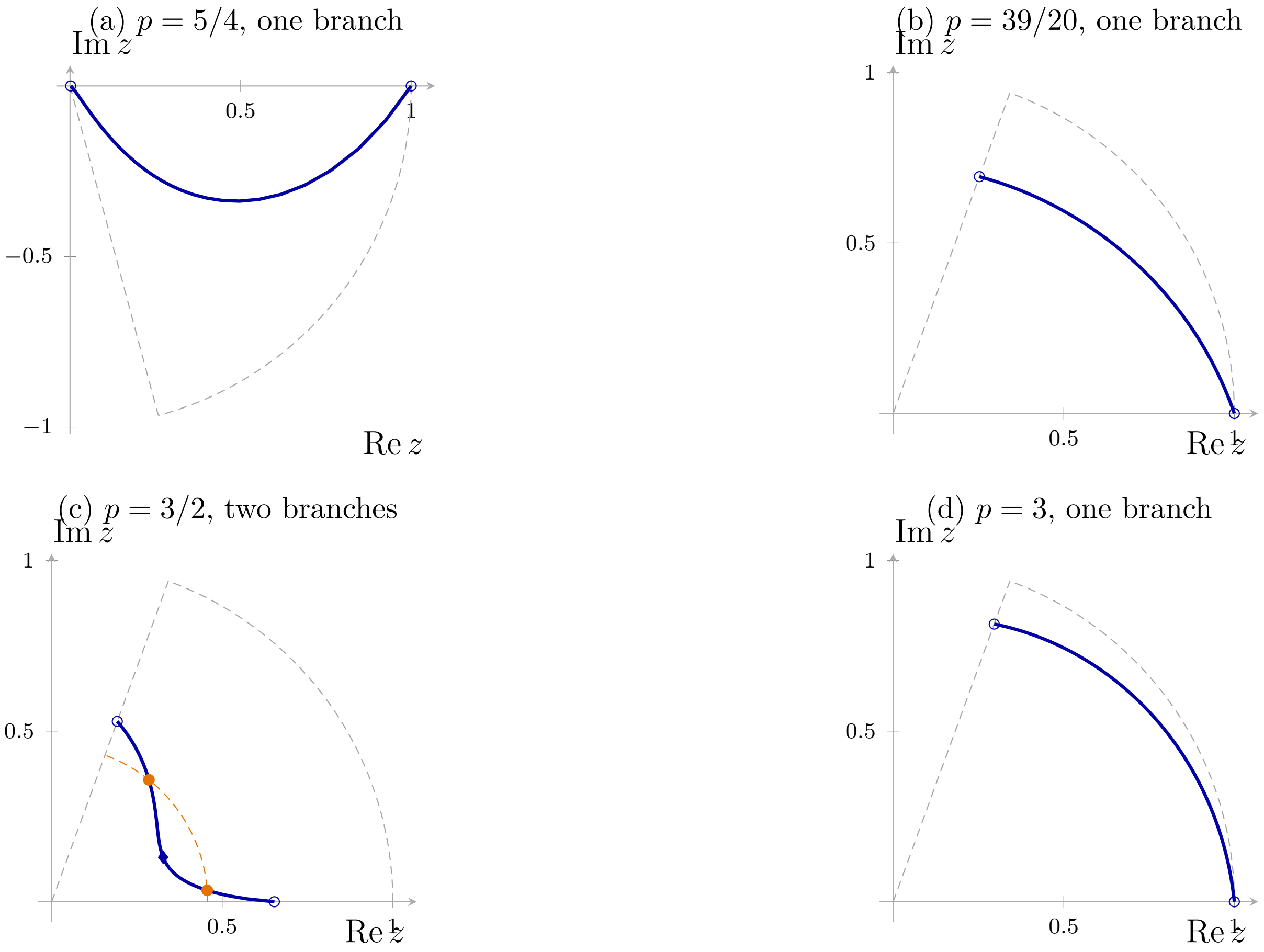
\begin{figure}[h!]
\centering
\begin{tikzpicture}
\begin{axis}[
width=.45\textwidth,height=.39\textwidth,
xmin=-.04,xmax=1.07,ymin=-1.02,ymax=.06,
axis lines=middle,axis line style={gray!65},
axis equal image,
xtick={.5,1},ytick={-.5,-1},
xlabel={$\operatorname{Re}z$},ylabel={$\operatorname{Im}z$},
xlabel style={at={(axis description cs:1,.04)},anchor=north east},
ylabel style={at={(axis description cs:.01,1)},anchor=south west},
tick label style={font=\scriptsize},
title={\small (a) $p=5/4$, one branch},
]
\addplot[gray!70,densely dashed,no marks] coordinates {(1.000000,0.000000) (0.999514,-0.031162) (0.998058,-0.062293) (0.995632,-0.093364) (0.992239,-0.124344) (0.987883,-0.155203) (0.982566,-0.185912) (0.976296,-0.216440) (0.969077,-0.246757) (0.960917,-0.276836) (0.951824,-0.306645) (0.941806,-0.336156) (0.930874,-0.365341) (0.919037,-0.394171) (0.906308,-0.422618) (0.892698,-0.450655) (0.878222,-0.478254) (0.862892,-0.505388) (0.846724,-0.532032) (0.829734,-0.558159) (0.811938,-0.583744) (0.793353,-0.608761) (0.773998,-0.633188) (0.753891,-0.656999) (0.733052,-0.680173) (0.711501,-0.702685) (0.689258,-0.724516) (0.666347,-0.745642) (0.642788,-0.766044) (0.618604,-0.785703) (0.593820,-0.804598) (0.568459,-0.822711) (0.542546,-0.840026) (0.516106,-0.856525) (0.489165,-0.872191) (0.461749,-0.887011) (0.433884,-0.900969) (0.405597,-0.914052) (0.376917,-0.926247) (0.347871,-0.937542) (0.318487,-0.947927) (0.288793,-0.957392) (0.258819,-0.965926)};
\addplot[gray!70,densely dashed,no marks] coordinates {(0.000000,0.000000) (0.258819,-0.965926)};
\addplot[blue!65!black,very thick,no marks] coordinates {(1.000000,0.000000) (0.924090,-0.102451) (0.844351,-0.185267) (0.764645,-0.247532) (0.688358,-0.290777) (0.617724,-0.317976) (0.553797,-0.332514) (0.496759,-0.337525) (0.446270,-0.335618) (0.401741,-0.328835) (0.362501,-0.318722) (0.327889,-0.306423) (0.297298,-0.292775) (0.270189,-0.278387) (0.246098,-0.263695) (0.224624,-0.249011) (0.205426,-0.234558) (0.188211,-0.220487) (0.172733,-0.206904) (0.158778,-0.193874) (0.146167,-0.181439) (0.134741,-0.169619) (0.124368,-0.158420) (0.114931,-0.147839) (0.106328,-0.137864) (0.098471,-0.128477) (0.091284,-0.119657) (0.084699,-0.111383) (0.078657,-0.103629) (0.073104,-0.096369) (0.067995,-0.089579) (0.063289,-0.083234) (0.058947,-0.077307) (0.054939,-0.071775) (0.051234,-0.066615) (0.047805,-0.061804) (0.044630,-0.057320) (0.041687,-0.053143) (0.038957,-0.049254) (0.036422,-0.045634) (0.034066,-0.042265) (0.031876,-0.039132) (0.029838,-0.036218) (0.027941,-0.033509) (0.026173,-0.030991) (0.024525,-0.028652) (0.022988,-0.026479) (0.021554,-0.024461) (0.020215,-0.022587) (0.018964,-0.020847) (0.017795,-0.019233) (0.016702,-0.017736) (0.015680,-0.016346) (0.014723,-0.015058) (0.013828,-0.013863) (0.012989,-0.012755) (0.012204,-0.011728) (0.011468,-0.010777) (0.010778,-0.009895) (0.010130,-0.009079) (0.009524,-0.008323) (0.008954,-0.007623) (0.008420,-0.006974) (0.007918,-0.006374) (0.007447,-0.005819) (0.007005,-0.005305) (0.006589,-0.004830) (0.006199,-0.004390) (0.005832,-0.003983) (0.005487,-0.003607) (0.005163,-0.003259) (0.004857,-0.002937) (0.004570,-0.002639) (0.004300,-0.002363) (0.004046,-0.002107) (0.003807,-0.001871) (0.003581,-0.001651) (0.003369,-0.001447) (0.003168,-0.001258) (0.002980,-0.001081) (0.002801,-0.000915) (0.002633,-0.000758) (0.002473,-0.000609) (0.002322,-0.000465) (0.002178,-0.000323) (0.002040,-0.000174) (0.001904,0.000000)};
\addplot[blue!65!black,only marks,mark=o,mark size=1.8pt] coordinates {(1.000000,0.000000) (0.001904,0.000000)};
\end{axis}
\end{tikzpicture}
\hfill
\begin{tikzpicture}
\begin{axis}[
width=.45\textwidth,height=.39\textwidth,
xmin=-.04,xmax=1.07,ymin=-.06,ymax=1.02,
axis lines=middle,axis line style={gray!65},
axis equal image,
xtick={.5,1},ytick={.5,1},
xlabel={$\operatorname{Re}z$},ylabel={$\operatorname{Im}z$},
xlabel style={at={(axis description cs:1,.04)},anchor=north east},
ylabel style={at={(axis description cs:.01,1)},anchor=south west},
tick label style={font=\scriptsize},
title={\small (b) $p=39/20$, one branch},
]
\addplot[gray!70,densely dashed,no marks] coordinates {(1.000000,0.000000) (0.999577,0.029085) (0.998308,0.058145) (0.996195,0.087156) (0.993238,0.116093) (0.989442,0.144932) (0.984808,0.173648) (0.979341,0.202218) (0.973045,0.230616) (0.965926,0.258819) (0.957990,0.286803) (0.949243,0.314545) (0.939693,0.342020) (0.929348,0.369206) (0.918216,0.396080) (0.906308,0.422618) (0.893633,0.448799) (0.880201,0.474600) (0.866025,0.500000) (0.851117,0.524977) (0.835488,0.549509) (0.819152,0.573576) (0.802123,0.597159) (0.784416,0.620235) (0.766044,0.642788) (0.747025,0.664796) (0.727374,0.686242) (0.707107,0.707107) (0.686242,0.727374) (0.664796,0.747025) (0.642788,0.766044) (0.620235,0.784416) (0.597159,0.802123) (0.573576,0.819152) (0.549509,0.835488) (0.524977,0.851117) (0.500000,0.866025) (0.474600,0.880201) (0.448799,0.893633) (0.422618,0.906308) (0.396080,0.918216) (0.369206,0.929348) (0.342020,0.939693)};
\addplot[gray!70,densely dashed,no marks] coordinates {(0.000000,0.000000) (0.342020,0.939693)};
\addplot[blue!65!black,very thick,no marks] coordinates {(1.000000,0.000000) (0.995006,0.014136) (0.989839,0.028131) (0.984502,0.041983) (0.979001,0.055691) (0.973338,0.069254) (0.967518,0.082669) (0.961545,0.095936) (0.955424,0.109053) (0.949158,0.122021) (0.942750,0.134837) (0.936206,0.147502) (0.929529,0.160014) (0.922723,0.172373) (0.915791,0.184578) (0.908737,0.196630) (0.901565,0.208528) (0.894278,0.220272) (0.886881,0.231861) (0.879376,0.243296) (0.871767,0.254577) (0.864056,0.265703) (0.856249,0.276676) (0.848347,0.287495) (0.840354,0.298161) (0.832272,0.308674) (0.824106,0.319034) (0.815857,0.329243) (0.807529,0.339300) (0.799124,0.349206) (0.790646,0.358962) (0.782096,0.368569) (0.773478,0.378028) (0.764793,0.387338) (0.756045,0.396502) (0.747235,0.405519) (0.738366,0.414391) (0.729441,0.423119) (0.720460,0.431703) (0.711427,0.440145) (0.702344,0.448445) (0.693212,0.456605) (0.684033,0.464625) (0.674809,0.472506) (0.665541,0.480250) (0.656233,0.487857) (0.646884,0.495329) (0.637497,0.502667) (0.628073,0.509871) (0.618614,0.516943) (0.609121,0.523884) (0.599595,0.530694) (0.590038,0.537375) (0.580451,0.543928) (0.570834,0.550354) (0.561190,0.556653) (0.551520,0.562828) (0.541823,0.568878) (0.532102,0.574805) (0.522357,0.580610) (0.512589,0.586294) (0.502799,0.591857) (0.492988,0.597301) (0.483156,0.602626) (0.473305,0.607834) (0.463434,0.612925) (0.453546,0.617900) (0.443639,0.622761) (0.433715,0.627507) (0.423775,0.632139) (0.413819,0.636659) (0.403847,0.641067) (0.393859,0.645363) (0.383858,0.649549) (0.373842,0.653625) (0.363812,0.657592) (0.353768,0.661450) (0.343712,0.665201) (0.333643,0.668844) (0.323561,0.672379) (0.313467,0.675809) (0.303361,0.679133) (0.293243,0.682351) (0.283114,0.685464) (0.272974,0.688473) (0.262823,0.691378) (0.252660,0.694179)};
\addplot[blue!65!black,only marks,mark=o,mark size=1.8pt] coordinates {(1.000000,0.000000) (0.252660,0.694179)};
\end{axis}
\end{tikzpicture}

\medskip
\begin{tikzpicture}
\begin{axis}[
width=.45\textwidth,height=.39\textwidth,
xmin=-.04,xmax=1.07,ymin=-.06,ymax=1.02,
axis lines=middle,axis line style={gray!65},
axis equal image,
xtick={.5,1},ytick={.5,1},
xlabel={$\operatorname{Re}z$},ylabel={$\operatorname{Im}z$},
xlabel style={at={(axis description cs:1,.04)},anchor=north east},
ylabel style={at={(axis description cs:.01,1)},anchor=south west},
tick label style={font=\scriptsize},
title={\small (c) $p=3/2$, two branches},
]
\addplot[gray!70,densely dashed,no marks] coordinates {(1.000000,0.000000) (0.999577,0.029085) (0.998308,0.058145) (0.996195,0.087156) (0.993238,0.116093) (0.989442,0.144932) (0.984808,0.173648) (0.979341,0.202218) (0.973045,0.230616) (0.965926,0.258819) (0.957990,0.286803) (0.949243,0.314545) (0.939693,0.342020) (0.929348,0.369206) (0.918216,0.396080) (0.906308,0.422618) (0.893633,0.448799) (0.880201,0.474600) (0.866025,0.500000) (0.851117,0.524977) (0.835488,0.549509) (0.819152,0.573576) (0.802123,0.597159) (0.784416,0.620235) (0.766044,0.642788) (0.747025,0.664796) (0.727374,0.686242) (0.707107,0.707107) (0.686242,0.727374) (0.664796,0.747025) (0.642788,0.766044) (0.620235,0.784416) (0.597159,0.802123) (0.573576,0.819152) (0.549509,0.835488) (0.524977,0.851117) (0.500000,0.866025) (0.474600,0.880201) (0.448799,0.893633) (0.422618,0.906308) (0.396080,0.918216) (0.369206,0.929348) (0.342020,0.939693)};
\addplot[gray!70,densely dashed,no marks] coordinates {(0.000000,0.000000) (0.342020,0.939693)};
\addplot[blue!65!black,very thick,no marks] coordinates {(0.652846,0.000000) (0.574350,0.008160) (0.531337,0.015101) (0.501040,0.021367) (0.477650,0.027172) (0.458679,0.032635) (0.442810,0.037835) (0.429249,0.042827) (0.417480,0.047652) (0.407147,0.052342) (0.397991,0.056923) (0.389820,0.061417) (0.382485,0.065843) (0.375870,0.070216) (0.369881,0.074550) (0.364442,0.078857) (0.359491,0.083149) (0.354974,0.087434) (0.350846,0.091723) (0.347069,0.096023) (0.343610,0.100342) (0.340439,0.104687) (0.337532,0.109065) (0.334866,0.113482) (0.332420,0.117944) (0.330175,0.122456) (0.328116,0.127023) (0.326227,0.131650) (0.324494,0.136343) (0.322904,0.141105) (0.321445,0.145940) (0.320104,0.150852) (0.318872,0.155845) (0.317737,0.160922) (0.316690,0.166086) (0.315720,0.171339) (0.314819,0.176685) (0.313977,0.182125) (0.313185,0.187662) (0.312435,0.193297) (0.311716,0.199030) (0.311022,0.204864) (0.310343,0.210798) (0.309670,0.216833) (0.308996,0.222969) (0.308312,0.229206) (0.307610,0.235542) (0.306881,0.241976) (0.306118,0.248507) (0.305312,0.255133) (0.304456,0.261852) (0.303542,0.268661) (0.302562,0.275558) (0.301510,0.282538) (0.300377,0.289600) (0.299157,0.296738) (0.297842,0.303949) (0.296427,0.311228) (0.294905,0.318572) (0.293269,0.325975) (0.291515,0.333432) (0.289637,0.340938) (0.287628,0.348489) (0.285486,0.356078) (0.283204,0.363700) (0.280779,0.371350) (0.278206,0.379022) (0.275482,0.386710) (0.272605,0.394409) (0.269570,0.402113) (0.266375,0.409817) (0.263017,0.417514) (0.259496,0.425200) (0.255808,0.432868) (0.251952,0.440513) (0.247927,0.448130) (0.243732,0.455713) (0.239367,0.463257) (0.234830,0.470757) (0.230122,0.478208) (0.225243,0.485605) (0.220192,0.492942) (0.214970,0.500215) (0.209577,0.507420) (0.204015,0.514552) (0.198285,0.521605) (0.192386,0.528577)};
\addplot[blue!65!black,only marks,mark=o,mark size=1.8pt] coordinates {(0.652846,0.000000) (0.192386,0.528577)};
\addplot[orange!90!black,densely dashed,no marks] coordinates {(0.457140,0.000000) (0.456807,0.017449) (0.455808,0.034872) (0.454145,0.052245) (0.451820,0.069542) (0.448836,0.086737) (0.445199,0.103806) (0.440912,0.120723) (0.435983,0.137465) (0.430418,0.154006) (0.424226,0.170323) (0.417416,0.186391) (0.409997,0.202188) (0.401981,0.217690) (0.393378,0.232875) (0.384203,0.247721) (0.374468,0.262205) (0.364186,0.276307) (0.353374,0.290007) (0.342047,0.303284) (0.330222,0.316119) (0.317915,0.328493) (0.305145,0.340388) (0.291930,0.351787) (0.278289,0.362674) (0.264243,0.373032) (0.249812,0.382846) (0.235017,0.392102) (0.219879,0.400787) (0.204421,0.408888) (0.188665,0.416393) (0.172634,0.423291) (0.156351,0.429572)};
\addplot[orange!90!black,only marks,mark=*,mark size=1.9pt] coordinates {(0.455911,0.033498) (0.285120,0.357329)};
\addplot[blue!65!black,only marks,mark=diamond*,mark size=2.2pt] coordinates {(0.326731,0.130372)};
\end{axis}
\end{tikzpicture}
\hfill
\begin{tikzpicture}
\begin{axis}[
width=.45\textwidth,height=.39\textwidth,
xmin=-.04,xmax=1.07,ymin=-.06,ymax=1.02,
axis lines=middle,axis line style={gray!65},
axis equal image,
xtick={.5,1},ytick={.5,1},
xlabel={$\operatorname{Re}z$},ylabel={$\operatorname{Im}z$},
xlabel style={at={(axis description cs:1,.04)},anchor=north east},
ylabel style={at={(axis description cs:.01,1)},anchor=south west},
tick label style={font=\scriptsize},
title={\small (d) $p=3$, one branch},
]
\addplot[gray!70,densely dashed,no marks] coordinates {(1.000000,0.000000) (0.999577,0.029085) (0.998308,0.058145) (0.996195,0.087156) (0.993238,0.116093) (0.989442,0.144932) (0.984808,0.173648) (0.979341,0.202218) (0.973045,0.230616) (0.965926,0.258819) (0.957990,0.286803) (0.949243,0.314545) (0.939693,0.342020) (0.929348,0.369206) (0.918216,0.396080) (0.906308,0.422618) (0.893633,0.448799) (0.880201,0.474600) (0.866025,0.500000) (0.851117,0.524977) (0.835488,0.549509) (0.819152,0.573576) (0.802123,0.597159) (0.784416,0.620235) (0.766044,0.642788) (0.747025,0.664796) (0.727374,0.686242) (0.707107,0.707107) (0.686242,0.727374) (0.664796,0.747025) (0.642788,0.766044) (0.620235,0.784416) (0.597159,0.802123) (0.573576,0.819152) (0.549509,0.835488) (0.524977,0.851117) (0.500000,0.866025) (0.474600,0.880201) (0.448799,0.893633) (0.422618,0.906308) (0.396080,0.918216) (0.369206,0.929348) (0.342020,0.939693)};
\addplot[gray!70,densely dashed,no marks] coordinates {(0.000000,0.000000) (0.342020,0.939693)};
\addplot[blue!65!black,very thick,no marks] coordinates {(1.000000,0.000000) (0.998442,0.014185) (0.996684,0.028326) (0.994729,0.042420) (0.992576,0.056464) (0.990227,0.070455) (0.987682,0.084392) (0.984942,0.098270) (0.982009,0.112088) (0.978884,0.125842) (0.975567,0.139531) (0.972059,0.153150) (0.968362,0.166699) (0.964477,0.180173) (0.960405,0.193571) (0.956148,0.206889) (0.951706,0.220126) (0.947081,0.233278) (0.942274,0.246343) (0.937287,0.259318) (0.932121,0.272202) (0.926778,0.284991) (0.921259,0.297682) (0.915566,0.310275) (0.909699,0.322765) (0.903662,0.335151) (0.897455,0.347430) (0.891081,0.359599) (0.884540,0.371657) (0.877835,0.383602) (0.870968,0.395429) (0.863939,0.407139) (0.856752,0.418727) (0.849408,0.430193) (0.841909,0.441533) (0.834257,0.452745) (0.826454,0.463828) (0.818502,0.474780) (0.810403,0.485597) (0.802158,0.496278) (0.793772,0.506821) (0.785244,0.517224) (0.776578,0.527485) (0.767776,0.537602) (0.758839,0.547573) (0.749771,0.557396) (0.740574,0.567069) (0.731249,0.576591) (0.721800,0.585959) (0.712229,0.595172) (0.702537,0.604228) (0.692728,0.613125) (0.682804,0.621861) (0.672767,0.630436) (0.662621,0.638847) (0.652367,0.647092) (0.642008,0.655171) (0.631547,0.663082) (0.620986,0.670822) (0.610328,0.678392) (0.599576,0.685788) (0.588732,0.693011) (0.577800,0.700058) (0.566781,0.706929) (0.555679,0.713621) (0.544496,0.720135) (0.533235,0.726468) (0.521900,0.732620) (0.510492,0.738589) (0.499015,0.744374) (0.487472,0.749975) (0.475866,0.755390) (0.464199,0.760619) (0.452474,0.765660) (0.440695,0.770512) (0.428865,0.775176) (0.416986,0.779650) (0.405061,0.783932) (0.393094,0.788024) (0.381087,0.791923) (0.369045,0.795630) (0.356968,0.799144) (0.344862,0.802464) (0.332729,0.805590) (0.320572,0.808521) (0.308394,0.811257) (0.296198,0.813798)};
\addplot[blue!65!black,only marks,mark=o,mark size=1.8pt] coordinates {(1.000000,0.000000) (0.296198,0.813798)};
\end{axis}
\end{tikzpicture}
\caption{The blue curves are the coupling lines and the open circles mark the boundary points of the branches. In panel (c) the diamond marks the minimal radius, the orange dashed arc is the circle $r\simeq0.457$, and the two orange dots have the same radius but different angles.}
\label{fig:coupling-branches}
\end{figure}

Figure~\ref{fig:coupling-branches} shows the coupling lines in the input plane $z=re^{i\theta}$. The dashed curves mark the unit circle and the corresponding sector boundaries. Panel (a) uses $\rho=1/2$ and $\gamma=75^\circ$, and the other three panels use $\rho=3/4$ and $\gamma=70^\circ$. Panel (a) is computed by solving for the angle at given radius, the other panels by solving for the radius at given angle, in each case by bisection on \eqref{eq:coupling-hyp}.

In panel (a) the branch in the fourth quadrant is single-valued on $(e^{-2\tau_*/p},1)$. In panels (b) and (d) the radius in the first quadrant strictly decreases with the angle. Panel (c) satisfies \eqref{eq:coupling-fold-criterion}: the radius first decreases and then increases, with minimal radius about $0.352$ and end radii about $0.653$ and $0.5625$. Hence a radius between the minimum and $0.5625$ corresponds to two angles. In particular, for $p>2$ the first quadrant carries a single branch only; the two-branch phenomenon in the first quadrant appears for $1<p<2$.

\clearpage
\section*{Acknowledgements}

The author acknowledges the use of AI-assisted tools in the preparation of this manuscript.
All mathematical results have been verified by the author,
who takes full responsibility for the content of the paper.

\end{document}